\documentclass[a4paper,11pt,reqno]{amsart}
\usepackage{amsmath,amsthm,amssymb,amscd}
\usepackage[margin=2.3cm]{geometry}
\usepackage{xcolor}
\usepackage{parskip}
\usepackage{mathtools}
\usepackage{booktabs}
\usepackage{placeins}
\usepackage{thmtools}
\mathtoolsset{showonlyrefs}
\usepackage{csquotes}
\usepackage[activate={true,nocompatibility},final,tracking=true,kerning=true,spacing=true]{microtype}
\usepackage[pdfencoding=unicode,pdftex,bookmarks=true,bookmarksnumbered]{hyperref}

\usepackage[shortlabels]{enumitem}
\setlist[enumerate,1]{label={\upshape(\roman*)}}

\definecolor{citecolour}{rgb}{0.0,0.0,0.8}
\colorlet{linkcolour}{green!50!black}
\hypersetup{
  colorlinks,
  breaklinks,
  linkcolor=linkcolour,
  citecolor=citecolour,
  filecolor=linkcolour,
  urlcolor=linkcolour
}

\theoremstyle{plain}
\newtheorem{theorem}{Theorem}[section]
\newtheorem{proposition}[theorem]{Proposition}
\newtheorem{lemma}[theorem]{Lemma}
\newtheorem{corollary}[theorem]{Corollary}
\newtheorem{question}[theorem]{Question}
\newtheorem{theoremA}{Theorem}

\theoremstyle{remark}
\newtheorem{remark}[theorem]{Remark}
\newtheorem{example}[theorem]{Example}

\numberwithin{equation}{section}

\newcommand{\card}[1]{\left\lvert #1\right\rvert}
\newcommand{\gauss}[2]{\genfrac{[}{]}{0pt}{}{#1}{#2}_p}
\newcommand{\qgauss}[2]{\genfrac{[}{]}{0pt}{}{#1}{#2}_q}
\renewcommand{\leq}{\leqslant}
\renewcommand{\geq}{\geqslant}
\renewcommand{\emptyset}{\varnothing}

\newenvironment{proofof}{{\bf Proof. }}{\hfill $\blacksquare$ \\}

\begin{document}

\title{Subgroup bounds for abelian $p$-groups with applications}
\author{Stefanos Aivazidis}
\email{stefanosaivazidis@gmail.com}
\date{\today}

\begin{abstract}
We prove an upper bound for the number of subgroups of an abelian $p$-group,
with constants that are sharp on homocyclic blocks. If $A_\lambda$ has type
$\lambda$, the main exponent is
$\eta(\lambda)=\sum_i\lfloor(\lambda'_i)^2/4\rfloor$, where the parts of
$\lambda'$ are the column heights of the Ferrers diagram of $\lambda$. 
The leading coefficient is the number of central choices in the blocks of
columns of odd height. 
The proof starts from Delsarte's formula, separates the blocks by column
height, and estimates each homocyclic block by successive Durfee squares.
We also give three applications: a comparison theorem for $p$-groups of
nilpotency class less than $p$, diagonal summatory estimates with the shape
fixed, and formulae for leading terms at fixed rank for finite abelian
$p$-groups.
\end{abstract}

\maketitle

\section{Introduction}
\label{sec:overview}

The number of subgroups of a finite $p$-group is one of the most natural
arithmetic invariants one can attach to it, and one of the more stubborn to
describe. The point of departure of the work reported here was our desire
to have, for a finite $p$-group, an upper bound for this number that responds
to some of the internal structure of the group and not merely to its order. For groups whose
nilpotency class is less than $p$ such a bound is within reach; the route
to it, due to Groves, passes through an associated abelian $p$-group and so
throws the whole weight of the problem onto the abelian case. It is the
abelian case that we settle here.

Let $p$ be a prime, and let $A$ be a finite abelian $p$-group; we write
$s(A)$ for the number of subgroups of $A$. By the classification theorem,
$A$ is uniquely of the form
\[
        A_\lambda=
        C_{p^{\lambda_1}}\times C_{p^{\lambda_2}}\times\cdots,
\]
where $\lambda=(\lambda_1,\lambda_2,\ldots)$ is a partition. Formulae for
the number of subgroups of each given type go back to
Delsarte~\cite{Delsarte} and Birkhoff~\cite{Birkhoff}; exact as they are,
they do not by themselves display a compact global bound for the total
number of subgroups. What we isolate here is a bound of that kind, sensitive
to the shape of $\lambda$, with an explicit leading term and explicit block
constants.

The shape relevant to the estimate is most naturally expressed through the
conjugate partition $\lambda'$, whose parts are the column heights of the
Ferrers diagram of $\lambda$. 
In Delsarte's formula, recalled in \eqref{eq:delsarte}, 
the leading exponent separates over these columns: 
a column of height $h$ in which the subgroup column has height $x$ 
contributes $x(h-x)$ to the leading power of $p$. 
The largest contribution of this column is therefore $\lfloor h^2/4\rfloor$, 
and the largest exponent occurring in the total subgroup count is
\[
        \eta(\lambda) \coloneqq
        \sum_{i\geq 1}\left\lfloor\frac{(\lambda'_i)^2}{4}\right\rfloor.
\]
For each height $h\geq 1$, put $d_h(\lambda)\coloneqq\lambda_h-\lambda_{h+1}$, with
trailing zeros appended; then $d_h(\lambda)$ is the number of columns of
height exactly $h$. Set
\[
        \theta(\lambda) \coloneqq
        \prod_{\substack{h\geq 1\\ h\textnormal{ odd}}}(d_h(\lambda)+1),
        \quad \text{and} \quad
        r(\lambda) \coloneqq \card{\{h\geq 2:d_h(\lambda)>0\}}.
\]
A column of even height has a single value of $x$ realising
$x(h-x)=\lfloor h^2/4\rfloor$, while a column of odd height has two; in a
block of $d_h(\lambda)$ equal columns of odd height, monotonicity admits a
single switch between the two values, giving $d_h(\lambda)+1$ choices. 
The factor $\theta(\lambda)$ counts these maximising types, and we prove in
Section~\ref{sec:delsarte} that $\theta(\lambda)p^{\eta(\lambda)}$ is the
exact leading term of $s(A_\lambda)$ as a polynomial in $p$.

Put $q\coloneqq p^{-1}$, and set $(q;q)_a\coloneqq\prod_{j=1}^a(1-q^j)$, with
$(q;q)_0\coloneqq 1$. The block constants in the upper bound below are
\begin{equation}\label{eq:block-constants}
\begin{aligned}
        \xi_{2a}(q)
        &=
        \frac{(q;q)_{2a}}{(q;q)_a^4}
        &&(a\geq 1), \\
        \xi_{2a+1}(q)
        &=
        (1-q)\frac{(q;q)_{2a+1}}{(q;q)_a^2(q;q)_{a+1}^2}
        &&(a\geq 0).
\end{aligned}
\end{equation}
The second formula gives $\xi_1(q)=1$, so a block of height one contributes no
nontrivial factor. Thus the product in Theorem~\ref{thm:main} may be taken
from $h\geq 2$ without loss.

Our main result is the following.

\begin{theoremA}\label{thm:main}
For every partition $\lambda$,
\[
        \theta(\lambda)p^{\eta(\lambda)}
        \leq
        s(A_\lambda)
        \leq
        \theta(\lambda)p^{\eta(\lambda)}
        \prod_{\substack{h\geq 2\\ d_h(\lambda)>0}}\xi_h(q).
\]
Moreover, for fixed $\lambda$,
\[
        s(A_\lambda)=
        \theta(\lambda)p^{\eta(\lambda)}
        +O_\lambda(p^{\eta(\lambda)-1})
        \quad \textnormal{as } p\to\infty.
\]
\end{theoremA}

Thus the lower bound is the exact leading term of $s(A_\lambda)$ as a
polynomial in $p$, and the second assertion records the corresponding
fixed-$\lambda$ asymptotic. The upper bound
multiplies the leading term by one explicit factor $\xi_h(q)$ for each
nontrivial height block of height at least $2$. Bounding each $\xi_h(q)$
uniformly by $(q;q)_\infty^{-3}$ gives the coarser inequality
\begin{equation}\label{Eq:uniform-bound}
        s(A_\lambda)
        \leq
        K_p^{3r(\lambda)}\theta(\lambda)p^{\eta(\lambda)},
        \quad \text{where} \quad
        K_p\coloneqq(q;q)_\infty^{-1},
\end{equation}
which we shall use whenever the individual heights are immaterial.

A small example is worth recording before turning to the structure of the proof.
If $\lambda=(5,5,3,1)$, so that
$A_\lambda=C_{p^5}\times C_{p^5}\times C_{p^3}\times C_p$, then
$\lambda'=(4,3,3,2,2)$ and $\eta(\lambda)=10$.
The nonzero height blocks 
have widths $d_2(\lambda)=2$, $d_3(\lambda)=2$, and $d_4(\lambda)=1$, 
so $\theta(\lambda)=3$ and $r(\lambda)=3$.
Theorem~\ref{thm:main} therefore identifies $3p^{10}$ as the leading term
of $s(A_\lambda)$, and \eqref{Eq:uniform-bound} gives the uniform bound
$3K_p^9p^{10}$. The explicit block constants $\xi_2,\xi_3,\xi_4$ and the
complete polynomial $s(A_\lambda)$ for this $\lambda$ are computed in
Section~\ref{sec:delsarte}.

The proof of Theorem~\ref{thm:main} proceeds in three steps. First, once
Delsarte's formula is rewritten with $q=p^{-1}$, the leading power of $p$
in each summand separates from a bounded factor in $q$; the largest power
that occurs is $p^{\eta(\lambda)}$, attained by exactly $\theta(\lambda)$
summands, so the lower bound and the asymptotic for fixed $\lambda$ are immediate.
The whole difficulty lies in the upper bound, that is, in controlling the
normalised sum that remains. Here everything turns on a single observation:
a Ferrers diagram is not the product of its blocks of constant height, yet
the monotonicity conditions that tie adjacent blocks together may, once the
sum has been normalised, simply be dropped, since they only constrain a sum
of nonnegative terms. The problem is thereby reduced to a single homocyclic block. 
The estimate for such a block is obtained by decomposing Ferrers diagrams 
into successive Durfee squares; this yields the constants $\xi_h(q)$. 
These three steps, namely reading off the leading term,
estimating one homocyclic block, and decoupling the blocks of a general
diagram, occupy Sections~\ref{sec:delsarte}--\ref{sec:block-decoupling}.

The sections that follow give three applications of Theorem~\ref{thm:main}.
Section~\ref{sec:class-less-than-p} treats groups of nilpotency class less
than $p$, where Groves's comparison theorem reduces subgroup counting to the
abelian group determined by the logarithmic $\Omega$-profile;
Theorem~\ref{thm:main} then gives the corresponding bound.
Section~\ref{sec:diagonal-sums} fixes the shape $\lambda$ and lets the base
integer vary; the result is a multiplicative arithmetic function whose diagonal 
summatory function is governed by $\eta(\lambda)$ and $\theta(\lambda)$. 
Section~\ref{sec:fixed-rank} returns to a single prime and reads off leading
coefficients for abelian $p$-groups of fixed rank.
We close, in Section~\ref{sec:open-problems}, by recording three problems 
that we have been unable to settle.

\subsection{Notation and conventions}
\label{subsec:notation}

We collect here the notation and conventions used throughout the paper. Unless
otherwise stated, $p$ denotes a prime and all groups are finite. 
For a finite group $G$, write $s(G)$ for the number of subgroups of $G$, 
and write $C_m$ for the cyclic group of order $m$.

All partitions are written in nonincreasing order, and we freely adjoin trailing
zeros to partitions and to their conjugates. Thus expressions such as
$\lambda_i$, $\lambda'_i$, and $\mu'_i$ are defined for all $i\geq 1$.
Sums and products indexed by $i\geq 1$ and involving partition data are
understood with this convention of extending by zeros; in particular, all but
finitely many terms are trivial. In partition notation, $(m^r)$ denotes
$r$ parts equal to $m$, and terms with exponent $0$ are omitted. This is
separate from group notation such as $C_{p^d}^{\,h}$, where the superscript
denotes a direct product of $h$ copies.

For group-theoretic notation we use the following conventions.
The commutator is $[x,y]=xyx^{-1}y^{-1}$. When an action is needed, it is
specified explicitly. We write $Z(G)$, $G'$, $\Phi(G)$, and $\gamma_i(G)$
for the centre, derived subgroup, Frattini subgroup, and lower central series
of $G$, respectively.

If $G$ is a $p$-group, then
$\Omega_i(G)$ denotes the subgroup generated by the elements of $G$ whose
orders divide $p^i$. 
If $\exp(G)=p^e$, put
$\omega_i(G)=\log_p|\Omega_i(G)|$ for $0\leq i\leq e$, with
$\omega_0(G)=0$, and write
\[
\omega(G)=(\omega_0(G),\omega_1(G),\ldots,\omega_e(G))\,;
\]
we call $\omega(G)$ the (logarithmic) $\Omega$-profile of $G$.

We use standard asymptotic notation. In particular, we write $f\ll g$ or
$f=O(g)$ to mean $|f|\leq Cg$ for some constant $C$. A subscript on either
symbol indicates parameters on which $C$ is allowed to depend.

The main recurring notation is collected in
Tables~\ref{tab:notation} and \ref{tab:q-notation}. Local notation not listed there is defined
where it is used.

\FloatBarrier
\begin{table}[htbp]
\centering
\footnotesize
\begingroup
\renewcommand{\arraystretch}{1.12}
\begin{tabular}{@{}p{0.27\textwidth}p{0.66\textwidth}@{}}
\toprule
Notation & Meaning \\
\midrule
$\lambda=(\lambda_1,\lambda_2,\ldots)$
& Partition, nonincreasing; trailing zeros may be adjoined \\

$\ell(\lambda)$
& The number of nonzero parts of $\lambda$ \\

$|\lambda|$
& The size of $\lambda$, namely $\sum_i\lambda_i=\sum_i\lambda'_i$ \\

$\lambda'$
& Conjugate partition: $\lambda'_i=\card{\{j:\lambda_j\geq i\}}$ \\

$(m^r)$
& Partition notation for $r$ parts equal to $m$; terms with exponent $0$
  are omitted \\

$b\lambda$
& The partition $(b\lambda_1,b\lambda_2,\ldots)$; $0\lambda$ is the empty
  partition \\

$\mu\subseteq\lambda$
& Inclusion of Ferrers diagrams; equivalently, $\mu_i\leq\lambda_i$ for all
  $i$ \\

$A_\lambda$
& The abelian $p$-group
  $C_{p^{\lambda_1}}\times C_{p^{\lambda_2}}\times\cdots$ \\

$\alpha_\lambda(\mu;p)$
& The number of subgroups of $A_\lambda$ of type $\mu$ \\

$E_\lambda(\mu)$
& The Delsarte leading exponent
  $E_\lambda(\mu)=\sum_i\mu'_i(\lambda'_i-\mu'_i)$ \\

$d_h(\lambda)$
& Columns of height exactly $h$: $d_h(\lambda)=\lambda_h-\lambda_{h+1}$ \\

$\eta(\lambda)$
& The main exponent
  $\eta(\lambda)=\sum_i\lfloor(\lambda'_i)^2/4\rfloor$ \\

$\theta(\lambda)$
& The central-pattern factor
  $\prod_{h\geq 1,\ h\textnormal{ odd}}(d_h(\lambda)+1)$ \\

$r(\lambda)$
& Nonzero height blocks of height $\geq 2$:
  $\card{\{h\geq 2:d_h(\lambda)>0\}}$ \\
\bottomrule
\end{tabular}
\endgroup
\vspace{0.2em}
\caption{Partition and abelian $p$-group notation.}
\label{tab:notation}
\end{table}

\FloatBarrier
\begin{table}[htbp]
\centering
\footnotesize
\begingroup
\renewcommand{\arraystretch}{1.12}
\begin{tabular}{@{}p{0.27\textwidth}p{0.66\textwidth}@{}}
\toprule
Notation & Meaning \\
\midrule
$q$
& The reciprocal $p^{-1}$ \\

$(q;q)_a$, $(q;q)_\infty$
& The finite product $\prod_{j=1}^a(1-q^j)$, with $(q;q)_0=1$, and its
  infinite analogue \\

$\gauss{n}{m}$
& The Gaussian coefficient over $\mathbb F_p$ \\

$\qgauss{n}{m}$
& The $q$-Gaussian coefficient
  $(q;q)_n/((q;q)_m(q;q)_{n-m})$ \\

$K_p$
& The infinite product $(q;q)_\infty^{-1}$ \\

$\xi_h(q)$
& The height-$h$ homocyclic block constant in \eqref{eq:block-constants} \\

$\mathbf y=(y_1,\ldots,y_t)$
& A finite nonincreasing chain used in the Durfee-square estimates \\

$\mathcal D_a(\mathbf y)$
& The Durfee-chain product
  $\prod_{j=1}^t\qgauss{y_{j-1}}{y_j}$, with $y_0=a$ \\

$\Pi_h(x_1,\ldots,x_d)$
& The homocyclic block product
  $\prod_{j=1}^{d}\qgauss{h-x_{j+1}}{x_j-x_{j+1}}$, with $x_{d+1}=0$ \\
\bottomrule
\end{tabular}
\endgroup
\vspace{0.2em}
\caption{$q$-series, Durfee-chain, and homocyclic-block notation.}
\label{tab:q-notation}
\end{table}

\section{Delsarte's formula and the central terms}
\label{sec:delsarte}

We begin the proof of Theorem~\ref{thm:main} by recording Delsarte's
enumeration formula in the form needed here. If $\mu\subseteq\lambda$, let
$\alpha_\lambda(\mu;p)$ be the number of subgroups of $A_\lambda$ isomorphic
to $A_\mu$. The formula below expresses this number in terms of the conjugate
partitions $\lambda'$ and $\mu'$. It goes back to Delsarte~\cite{Delsarte};
in the present notation, see Butler~\cite[p.~21, Eq.~(1.9)]{Butler}. It may
also be viewed as Birkhoff's subgroup enumeration formula refined by
type~\cite{Birkhoff}.

\begin{proposition}\label{prop:delsarte}
For $\mu\subseteq\lambda$,
\begin{equation}\label{eq:delsarte}
        \alpha_\lambda(\mu;p)=
        \prod_{i\geq 1}
        p^{\mu'_{i+1}(\lambda'_i-\mu'_i)}
        \gauss{\lambda'_i-\mu'_{i+1}}{\mu'_i-\mu'_{i+1}}.    
\end{equation}
\end{proposition}

To extract the leading power of $p$ from each summand we use the factorisation
\begin{equation}\label{eq:gauss-factor}
        \gauss{n}{m}=p^{m(n-m)}\qgauss{n}{m},
\end{equation}
which splits each $\gauss{\cdot}{\cdot}$ in Delsarte's formula into a leading
power of $p$ and a $q$-Gaussian factor that is bounded as $p\to\infty$.
Applied to the $i$-th factor with $n=\lambda'_i-\mu'_{i+1}$ and
$m=\mu'_i-\mu'_{i+1}$, the relation \eqref{eq:gauss-factor} contributes
\[
        \mu'_{i+1}(\lambda'_i-\mu'_i)
        +(\mu'_i-\mu'_{i+1})(\lambda'_i-\mu'_i)
        =\mu'_i(\lambda'_i-\mu'_i)
\]
to the exponent of $p$. The leading exponent in $\alpha_\lambda(\mu;p)$ is
therefore
\[
        E_\lambda(\mu) \coloneqq
        \sum_i \mu'_i(\lambda'_i-\mu'_i).
\]
For a single column, the function $x(h-x)$ is maximised when $x$ is as close
as possible to $h/2$, whence
\[
        E_\lambda(\mu)
        \leq
        \sum_i\left\lfloor\frac{(\lambda'_i)^2}{4}\right\rfloor
        =\eta(\lambda).
\]

The condition $E_\lambda(\mu)=\eta(\lambda)$ is columnwise. For a column of
height $h$, the function $x(h-x)$ is maximised at $x=h/2$ if $h$ is even, and
at the two central values if $h$ is odd. Thus a block of height $h=2a$ forces
all corresponding entries of $\mu'$ to be $a$. A block of height $h=2a+1$ and
width $d_h(\lambda)$ admits exactly the monotone central patterns
\[
        \underbrace{a+1,\ldots,a+1}_{\text{$k$ entries}},
        \underbrace{a,\ldots,a}_{d_h(\lambda)-k\text{ entries}},
        \qquad 0\leq k\leq d_h(\lambda),
\]
and hence contributes $d_h(\lambda)+1$ choices. These choices are independent
from block to block: if $h>h'$, then
$\lfloor h/2\rfloor\geq \lceil h'/2\rceil$, so the patterns from consecutive
height blocks concatenate to a nonincreasing sequence. Therefore
\[
        \card{\{\mu\subseteq\lambda : E_\lambda(\mu)=\eta(\lambda)\}}
        =
        \prod_{\substack{h\geq 1\\ h\textnormal{ odd}}}
        (d_h(\lambda)+1)
        =
        \theta(\lambda).
\]

The Gaussian coefficient $\gauss{n}{m}$ is a monic polynomial in $p$ of
degree $m(n-m)$ with nonnegative coefficients; this is immediate from
\eqref{eq:gauss-factor}, since the $q$-binomial coefficient is a polynomial
in $q$ with nonnegative coefficients and constant term~$1$
(Andrews--Eriksson~\cite[Chap.~7]{AndrewsEriksson}). Consequently every
Delsarte summand with $E_\lambda(\mu)=\eta(\lambda)$ contributes a leading
term $p^{\eta(\lambda)}$, while every other summand has degree at most
$\eta(\lambda)-1$. Summing over the $\theta(\lambda)$ central types and
absorbing the finitely many noncentral summands into the error term, we
obtain
\[
        s(A_\lambda)=
        \theta(\lambda)p^{\eta(\lambda)}
        +O_\lambda(p^{\eta(\lambda)-1}).
\]
This establishes the lower bound and the fixed-$\lambda$ asymptotic of
Theorem~\ref{thm:main}.

We return to the example $\lambda=(5,5,3,1)$ from
Section~\ref{sec:overview}, with $\lambda'=(4,3,3,2,2)$, $\eta(\lambda)=10$,
and $\theta(\lambda)=3$. The only odd-height block contributing to
$\theta(\lambda)$ has height~$3$ and width~$2$; its three central choices
give the types
\[
        \mu'=(2,2,2,1,1),\qquad
        \mu'=(2,2,1,1,1),\qquad
        \mu'=(2,1,1,1,1),
\]
each with $E_\lambda(\mu)=10$, while every other $\mu\subseteq\lambda$
satisfies $E_\lambda(\mu)\leq 9$. Summing Delsarte's formula over all
$\mu\subseteq\lambda$ gives the exact polynomial
\[
\begin{aligned}
        s(A_\lambda)
        ={ }&3p^{10}+14p^9+30p^8+37p^7+40p^6+41p^5 \\
        &\quad +38p^4+33p^3+24p^2+13p+15.
\end{aligned}
\]
The leading term $3p^{10}$ confirms Theorem~\ref{thm:main}, and the
coefficient $14$ of $p^9$ shows that the error term
$O_\lambda(p^{\eta(\lambda)-1})$ cannot in general be improved.

The upper bound in Theorem~\ref{thm:main} requires the homocyclic
block constants for the heights occurring in $\lambda'$, namely
$h=2,3,4$. From \eqref{eq:block-constants},
\[
        \xi_2(q)
        = \frac{(q;q)_2}{(q;q)_1^4}
        = \frac{1-q^2}{(1-q)^3},
        \qquad
        \xi_3(q)
        = (1-q)\frac{(q;q)_3}{(q;q)_1^2(q;q)_2^2}
        = \frac{1-q^3}{(1-q)^2(1-q^2)},
\]
and
\[
        \xi_4(q)
        = \frac{(q;q)_4}{(q;q)_2^4}
        = \frac{(1-q^3)(1-q^4)}{(1-q)^3(1-q^2)^3}.
\]
Multiplying, Theorem~\ref{thm:main} gives
\[
        3p^{10}
        \leq s(A_\lambda)
        \leq
        3p^{10}\,\xi_2(q)\xi_3(q)\xi_4(q)
        =
        3p^{10}\,
        \frac{(1-q^3)^2(1-q^4)}{(1-q)^8(1-q^2)^3},
\]
while the uniform bound \eqref{Eq:uniform-bound} yields the cruder but more compact estimate
$s(A_\lambda)\leq 3K_p^9p^{10}$.

We shall use the following normalised form of a Delsarte summand in the sequel.
Combining Proposition~\ref{prop:delsarte} with \eqref{eq:gauss-factor} gives
\begin{equation}\label{eq:normalised-delsarte-summand}
        p^{-\eta(\lambda)}\alpha_\lambda(\mu;p)
        =
        q^{\sum_{i\geq 1}(\lfloor(\lambda'_i)^2/4\rfloor
        -\mu'_i(\lambda'_i-\mu'_i))}
        \prod_{i\geq 1}
        \qgauss{\lambda'_i-\mu'_{i+1}}{\mu'_i-\mu'_{i+1}}.
\end{equation}

\section{Successive Durfee squares}
\label{sec:durfee}

The upper bound will be proved by cutting the columns of the Ferrers diagram of
$\lambda$ into blocks of equal height. 
On each such block, the corresponding entries of the conjugate subgroup type
$\mu'$ form a finite nonincreasing sequence. After the middle-level insertion
used in Section~\ref{sec:homocyclic}, the sums to be estimated are expressed
in terms of such sequences. We shall call
them chains. This section records two estimates for the resulting chain sums.
The first is an exact generating-function identity supplied by successive
Durfee squares. The second is a shifted variant needed for blocks of odd height:
in that case entries equal to $1$ have zero shifted weight, and terminal strings
of ones account for the factor $N+1$.

We first fix the notation common to both estimates. For a finite nonincreasing
chain $\mathbf y=(y_1,\ldots,y_t)$ with
$a\geq y_1\geq\cdots\geq y_t\geq 1$, put $y_0=a$ and
\[
        \mathcal D_a(\mathbf y) \coloneqq
        \prod_{j=1}^t\qgauss{y_{j-1}}{y_j}.
\]
The empty chain is allowed; for it, $t=0$ and $\mathcal D_a(\emptyset)=1$. In
all sums over $y_1,\ldots,y_t$ below, $\mathbf y$ denotes the chain
$(y_1,\ldots,y_t)$, and the case $t=0$ is included.

\begin{lemma}\label{lem:durfee-identities}
For every $a\geq 0$, and with $z$ an auxiliary variable,
\begin{equation}\label{eq:durfee-chain-generating-function}
\sum_{\substack{t\geq 0\\ a\geq y_1\geq\cdots\geq y_t\geq 1}}
        z^{\sum_{j=1}^t y_j}q^{\sum_{j=1}^t y_j^2}
        \mathcal D_a(\mathbf y)
        =
        \prod_{r=1}^a\frac{1}{1-zq^r}.
\end{equation}
Consequently, the following identities hold.
\begin{enumerate}
\item
\begin{equation}\label{eq:durfee-chain-square-sum}
\sum_{\substack{t\geq 0\\ a\geq y_1\geq\cdots\geq y_t\geq 1}}
        q^{\sum_{j=1}^t y_j^2}
        \mathcal D_a(\mathbf y)
        =
        \frac{1}{(q;q)_a}.
\end{equation}

\item
\begin{equation}\label{eq:durfee-chain-shifted-square-sum}
\sum_{\substack{t\geq 0\\ a\geq y_1\geq\cdots\geq y_t\geq 1}}
        q^{\sum_{j=1}^t y_j(y_j+1)}
        \mathcal D_a(\mathbf y)
        =
        \prod_{r=2}^{a+1}\frac{1}{1-q^r}.
\end{equation}
\end{enumerate}
\end{lemma}

\begin{proofof}
The product on the right of
\eqref{eq:durfee-chain-generating-function} is the two-variable generating
function for partitions whose parts are at most $a$, with $q$ recording size
and $z$ recording the number of parts. Decompose such a partition into
successive Durfee squares, as in
Andrews--Eriksson~\cite[Chap.~8, Sec.~8.5]{AndrewsEriksson}: remove the
largest square in the upper-left corner, then repeat the construction on what
remains below it. If the successive side lengths are
$y_1,\ldots,y_t$, then the squares contribute $q^{\sum_j y_j^2}$ and their
rows contribute $z^{\sum_j y_j}$.

The portion of the Ferrers diagram to the right of the $j$-th square lies in a
rectangle with $y_j$ rows and $y_{j-1}-y_j$ columns. Its generating function
is therefore $\qgauss{y_{j-1}}{y_j}$. Multiplying over all successive squares
gives $\mathcal D_a(\mathbf y)$. Since the decomposition is unique and
reversible, including the empty chain which corresponds to the empty partition,
\eqref{eq:durfee-chain-generating-function} follows.

Putting $z=1$ in \eqref{eq:durfee-chain-generating-function} gives
\eqref{eq:durfee-chain-square-sum}. Putting $z=q$ changes the exponent
$\sum_j y_j^2$ into $\sum_j y_j(y_j+1)$, and the product becomes
\[
        \prod_{r=1}^a(1-q^{r+1})^{-1}
        =
        \prod_{r=2}^{a+1}(1-q^r)^{-1},
\]
which gives \eqref{eq:durfee-chain-shifted-square-sum}.
\end{proofof}

We shall also need the same type of estimate with the shifted exponent
$y(y-1)$. The difference is that entries equal to $1$ now have zero weight.
Thus terminal strings of ones may have any length up to $N$, and this is the
source of the factor $N+1$.

\begin{lemma}\label{lem:shifted-chain}
For $a\geq 1$ and every integer $N\geq 0$,
\[
\sum_{\substack{0\leq t\leq N\\ a\geq y_1\geq\cdots\geq y_t\geq 1}}
        q^{\sum_{j=1}^t y_j(y_j-1)}
        \mathcal D_a(\mathbf y)
        \leq
        \frac{N+1}{(q;q)_{a-1}}.
\]
Moreover,
\[
\lim_{N\to\infty}\frac{1}{N+1}
\sum_{\substack{0\leq t\leq N\\ a\geq y_1\geq\cdots\geq y_t\geq 1}}
        q^{\sum_{j=1}^t y_j(y_j-1)}
        \mathcal D_a(\mathbf y)
        =
        \frac{1}{(q;q)_{a-1}}.
\]
\end{lemma}

\begin{proofof}
Let $L_{a,N}^{(c)}$ be the part of the left side contributed by chains of
length at most $N$ with exactly $c$ terminal entries equal to $1$. Suppose
first that $c\geq 1$. Delete these terminal ones and subtract $1$ from every
remaining entry. What remains is a chain of length $s\leq N-c$,
\[
        a-1\geq z_1\geq\cdots\geq z_s\geq 1.
\]
The deleted terminal ones have zero shifted weight, since $1(1-1)=0$, and the
remaining shifted weight becomes
\[
        \sum_{i=1}^s (z_i+1)z_i=
        \sum_{i=1}^s z_i(z_i+1).
\]

We also need the Gaussian product after the same deletion. If $s\geq 1$, then
expanding the Gaussian coefficients in $q$-Pochhammer symbols gives
\[
\begin{aligned}
        \mathcal D_a(z_1+1,\ldots,z_s+1,\underbrace{1,\ldots,1}_{\text{$c$ entries}})
        &=
        \qgauss{a}{z_1+1}
        \prod_{i=2}^{s}\qgauss{z_{i-1}+1}{z_i+1}
        \qgauss{z_s+1}{1}                                      \\
        &=
        \frac{(q;q)_a}
        {(q;q)_1(q;q)_{a-1-z_1}
        \prod_{i=2}^{s}(q;q)_{z_{i-1}-z_i}(q;q)_{z_s}}           \\
        &=
        \qgauss{a}{1}\mathcal D_{a-1}(z_1,\ldots,z_s).
\end{aligned}
\]
The same identity holds when $s=0$, in which case both sides are
$\qgauss{a}{1}$. Thus, for fixed $c\geq 1$, the contribution is at most
\[
        \qgauss{a}{1}
        \sum_{\substack{t\geq 0\\ a-1\geq z_1\geq\cdots\geq z_t\geq 1}}
        q^{\sum_i z_i(z_i+1)}\mathcal D_{a-1}(z_1,\ldots,z_t),
\]
which, by \eqref{eq:durfee-chain-shifted-square-sum} with $a-1$ in place of $a$, is
\[
        \qgauss{a}{1}
        \prod_{r=2}^{a}\frac{1}{1-q^r}
        =
        \frac{1-q^a}{1-q}
        \prod_{r=2}^{a}\frac{1}{1-q^r}
        =
        \frac{1}{(q;q)_{a-1}}.
\]

It remains, for the inequality, to treat the chains with no terminal entry equal
to $1$. Appending one terminal $1$ to such a chain is an injective map into
the chains with exactly one terminal $1$, and the shifted exponent is
unchanged. If the original chain is empty, the Gaussian product is multiplied by
$\qgauss{a}{1}$; otherwise it is multiplied by $\qgauss{y_t}{1}$, where
$y_t\geq 2$ is the last entry of the original chain. In both cases the
multiplier is at least $1$, whence the $c=0$ contribution is bounded by the
unrestricted contribution of chains with one terminal $1$, which is bounded by
the quantity just computed. Since $0\leq c\leq N$, the asserted inequality
follows.

We turn to the limiting assertion. Put
\[
        R_M\coloneqq
        \sum_{\substack{0\leq s\leq M\\ a-1\geq z_1\geq\cdots\geq z_s\geq 1}}
        q^{\sum_i z_i(z_i+1)}\mathcal D_{a-1}(z_1,\ldots,z_s).
\]
By \eqref{eq:durfee-chain-shifted-square-sum} with $a-1$ in place of $a$,
the increasing sequence $R_M$ tends to
\[
        L=\prod_{r=2}^a(1-q^r)^{-1}.
\]
For the limiting assertion we keep the length restriction instead of discarding
it. Then, for $c\geq 1$, the deletion map above is exact, and therefore
$L_{a,N}^{(c)}=\qgauss{a}{1}R_{N-c}$. Hence
\[
        \sum_{c=1}^N L_{a,N}^{(c)}
        =
        \qgauss{a}{1}\sum_{M=0}^{N-1}R_M.
\]
The $c=0$ contribution is bounded independently of $N$, and so vanishes
after division by $N+1$. 
Since $R_M\to L$, the averages of the $R_M$'s also tend to $L$, and hence
\[
        \frac{1}{N+1}\sum_{M=0}^{N-1}R_M \longrightarrow L.
\]
Multiplying by $\qgauss{a}{1}$ gives
\[
        \qgauss{a}{1}L
        =
        \frac{1-q^a}{1-q}\prod_{r=2}^a\frac{1}{1-q^r}
        =
        \frac{1}{(q;q)_{a-1}},
\]
which is the stated limit.
\end{proofof}

\section{The homocyclic block}
\label{sec:homocyclic}

We are ready to estimate the homocyclic group $C_{p^d}^{\,h}$. A block of $d$ columns
of common height $h$ corresponds to this group, and the next section will cut a
general Ferrers diagram into such blocks. A subgroup type of $C_{p^d}^{\,h}$ is
encoded by a chain
\[
        h\geq x_1\geq\cdots\geq x_d\geq 0,
\]
where $x_j$ is the $j$-th part of the conjugate subgroup type. Set
$x_{d+1}=0$ and
\[
        \Pi_h(x_1,\ldots,x_d)\coloneqq
        \prod_{j=1}^{d}
        \qgauss{h-x_{j+1}}{x_j-x_{j+1}}.
\]
We shall use the following elementary insertion lemma. Inserting a level $m$
separates the entries above $m$ from those below it, and $\qgauss{h}{m}$
records the contribution of the inserted level.

\begin{lemma}\label{lem:middle-insertion}
Let $0\leq m\leq h$, and let
$\mathbf x=(x_1,\ldots,x_d)$ be a chain with
$h\geq x_1\geq\cdots\geq x_d\geq 0$. Form the upper chain
$\mathbf u=(u_1,\ldots,u_r)$ by taking, in order, the positive numbers
$u_i=x_i-m$ for which $x_i>m$. Form the lower chain
$\mathbf v=(v_1,\ldots,v_s)$ by taking the positive numbers $m-x_i$ for which
$x_i<m$, but in reverse order. Then
\[
        \Pi_h(\mathbf x)
        \leq
        \qgauss{h}{m}\mathcal D_{h-m}(\mathbf u)\mathcal D_m(\mathbf v).
\]
If one of the entries of $\mathbf x$ is equal to $m$, then equality holds.
\end{lemma}

\begin{proofof}
Writing the $q$-binomial coefficients in product form gives
\[
        \qgauss{h-x_{j+1}}{x_j-x_{j+1}}
        =
        \frac{(q;q)_{h-x_{j+1}}}
        {(q;q)_{x_j-x_{j+1}}(q;q)_{h-x_j}}.
\]
Multiplication over $j$ cancels the successive factors $(q;q)_{h-x_j}$, and
we obtain
\[
        \Pi_h(\mathbf x)=
        \frac{(q;q)_h}
        {(q;q)_{h-x_1}\prod_{j=1}^{d}(q;q)_{x_j-x_{j+1}}}.
\]
Suppose first that the value $m$ occurs in the chain. The drops from $h$ to
$m$ are precisely those recorded by the upper chain $\mathbf u$, while the drops
from $m$ to $0$ are those recorded by the lower chain $\mathbf v$, after the
order has been reversed, so the denominator above splits into the
denominators of $\qgauss{h}{m}$, $\mathcal D_{h-m}(\mathbf u)$, and
$\mathcal D_m(\mathbf v)$, and equality follows.

It remains to deal with the case where $m$ is absent. Insert $m$ into the
chain, allowing the insertion to occur at the top or at the bottom boundary.
Only the denominator factor corresponding to the jump across $m$ is affected.
Write this jump as a drop from $m+\rho$ to $m-\sigma$, where
$\rho,\sigma\geq 0$ and $\rho+\sigma>0$. The original factor is
$(q;q)_{\rho+\sigma}$, whereas after insertion it is replaced by
$(q;q)_\rho(q;q)_\sigma$. Since
\[
        (q;q)_{\rho+\sigma}
        =
        (q;q)_\rho\prod_{i=1}^{\sigma}(1-q^{\rho+i})
        \geq
        (q;q)_\rho\prod_{i=1}^{\sigma}(1-q^i)
        =
        (q;q)_\rho(q;q)_\sigma,
\]
the inserted product is not smaller than the original one. But the inserted
product is exactly
$\qgauss{h}{m}\mathcal D_{h-m}(\mathbf u)\mathcal D_m(\mathbf v)$.
\end{proofof}

\begin{proposition}\label{prop:homocyclic}
Let $d\geq 1$.  For $a\geq 1$,
\[
        s(C_{p^d}^{\,2a})
        \leq
        \xi_{2a}(q)p^{da^2}.
\]
For $a\geq 0$,
\[
        s(C_{p^d}^{\,2a+1})
        \leq
        (d+1)\xi_{2a+1}(q)p^{da(a+1)}.
\]
Both estimates are sharp as $d\to\infty$: for $a\geq 1$,
\[
        \lim_{d\to\infty}p^{-da^2}s(C_{p^d}^{\,2a})
        =
        \xi_{2a}(q),
\]
and for $a\geq 0$,
\[
        \lim_{d\to\infty}
        \frac{s(C_{p^d}^{\,2a+1})}{(d+1)p^{da(a+1)}}
        =
        \xi_{2a+1}(q).
\]
\end{proposition}

\begin{proofof}
The group $C_{p^d}^{\,h}$ has type $(d^h)$, whose conjugate partition has $d$
parts, each equal to $h$; thus $\eta\bigl((d^h)\bigr)=d\lfloor h^2/4\rfloor$.
Put
\[
        \delta_h(x)\coloneqq
        \left\lfloor\frac{h^2}{4}\right\rfloor-x(h-x).
\]
Specialising the normalised Delsarte summand of Section~\ref{sec:delsarte} to
this type and summing over all subgroup types $\mu\subseteq(d^h)$ (encoded,
as above, by the chains $h\geq x_1\geq\cdots\geq x_d\geq 0$) gives the
normalised identity
\[
        p^{-d\lfloor h^2/4\rfloor}s(C_{p^d}^{\,h})
        =
        \sum_{h\geq x_1\geq\cdots\geq x_d\geq 0}
        q^{\sum_{j=1}^{d}\delta_h(x_j)}\Pi_h(x_1,\ldots,x_d).
\]

Suppose first that $h=2a$. Then $\delta_{2a}(x)=(x-a)^2$. We insert the central
level $a$ by Lemma~\ref{lem:middle-insertion}. Entries above $a$ have the form
$a+u_i$, and entries below $a$ are encoded, in reverse order, by $a-v_i$; their
defects are $u_i^2$ and $v_i^2$, respectively. Thus a chain $\mathbf x$
determines upper and lower chains $\mathbf u=(u_1,\ldots,u_r)$ and
$\mathbf v=(v_1,\ldots,v_s)$ whose total length is at most $d$. Discarding this
length restriction only enlarges the sum, and therefore, with the empty chain
included in both sums,
\[
\begin{aligned}
        p^{-da^2}s(C_{p^d}^{\,2a})
        &\leq
        \qgauss{2a}{a}
        \sum_{\substack{r\geq 0\\ a\geq u_1\geq\cdots\geq u_r\geq 1}}
        q^{\sum_i u_i^2}\mathcal D_a(\mathbf u)  
        \sum_{\substack{s\geq 0\\ a\geq v_1\geq\cdots\geq v_s\geq 1}}
        q^{\sum_i v_i^2}\mathcal D_a(\mathbf v) \\
        &=
        \qgauss{2a}{a}\frac{1}{(q;q)_a^2}       \\
        &=
        \frac{(q;q)_{2a}}{(q;q)_a^4}            \\
        &=\xi_{2a}(q),
\end{aligned}
\]
where the two chain sums have been evaluated by
\eqref{eq:durfee-chain-square-sum}. This proves the even-height upper bound.

Suppose now that $h=2a+1$. If $a=0$, then $C_{p^d}^{\,h}$ is cyclic and
$s(C_{p^d})=d+1=(d+1)\xi_1(q)$. We may therefore assume that $a\geq 1$. Again
insert the level $a$. Let $\mathbf u=(u_1,\ldots,u_r)$ be the chain of positive
differences $x_i-a$, and let $\mathbf v=(v_1,\ldots,v_s)$ be the chain
obtained, in reverse order, from the positive differences $a-x_i$. Then
$1\leq u_i\leq a+1$, $1\leq v_i\leq a$, and a direct calculation gives
\[
        \delta_{2a+1}(a+u)=u(u-1),
        \qquad
        \delta_{2a+1}(a-v)=v(v+1).
\]
If the lower chain has length $s$, then the upper chain has length at most
$d-s$; entries equal to $a$, when present, fill the remaining positions and have
zero defect. Hence
\begin{multline*}
        p^{-da(a+1)}s(C_{p^d}^{\,2a+1})
        \leq
        \qgauss{2a+1}{a}
        \sum_{s=0}^{d}
        \sum_{a\geq v_1\geq\cdots\geq v_s\geq 1}
        q^{\sum_i v_i(v_i+1)}\mathcal D_a(\mathbf v)  \\
        {}\cdot
        \sum_{\substack{a+1\geq u_1\geq\cdots\geq u_r\geq 1\\0\leq r\leq d-s}}
        q^{\sum_i u_i(u_i-1)}\mathcal D_{a+1}(\mathbf u).
\end{multline*}

In view of Lemma~\ref{lem:shifted-chain}, 
the sum over the upper chain is at most $(d-s+1)/(q;q)_a$, 
and hence at most $(d+1)/(q;q)_a$. 
The remaining lower-chain sum is bounded by its unrestricted version, which by
\eqref{eq:durfee-chain-shifted-square-sum} is
\[
        \prod_{r=2}^{a+1}\frac{1}{1-q^r}
        =
        \frac{1-q}{(q;q)_{a+1}}.
\]
Consequently
\[
\begin{aligned}
        p^{-da(a+1)}s(C_{p^d}^{\,2a+1})
        &\leq
        \qgauss{2a+1}{a}
        \frac{d+1}{(q;q)_a}\frac{1-q}{(q;q)_{a+1}}  \\
        &=(d+1)(1-q)
        \frac{(q;q)_{2a+1}}{(q;q)_a^2(q;q)_{a+1}^2} \\
        &=(d+1)\xi_{2a+1}(q),
\end{aligned}
\]
which proves the odd-height upper bound.

It remains to prove sharpness. In the even case, restrict the sum to those
chains which contain at least one entry equal to $a$. Given upper and lower
chains of total length at most $d-1$, fill all remaining positions by entries
equal to $a$. On this restricted family Lemma~\ref{lem:middle-insertion} is an
equality. Any fixed pair of upper and lower chains is admitted for all
sufficiently large $d$, and all summands are nonnegative. Monotone convergence,
together with \eqref{eq:durfee-chain-square-sum}, gives
\[
\begin{aligned}
        \liminf_{d\to\infty}p^{-da^2}s(C_{p^d}^{\,2a})
        &\geq
        \qgauss{2a}{a}
        \left(\sum_{\substack{r\geq 0\\ a\geq u_1\geq\cdots\geq u_r\geq 1}}
        q^{\sum_i u_i^2}\mathcal D_a(\mathbf u)\right)^2  \\
        &=
        \qgauss{2a}{a}(q;q)_a^{-2}
        =\xi_{2a}(q).
\end{aligned}
\]
The upper bound gives the reverse inequality.

For the odd case, assume first that $a\geq 1$ and again restrict to chains
containing at least one entry equal to $a$. Let
\[
        U_N\coloneqq
        \sum_{\substack{a+1\geq u_1\geq\cdots\geq u_r\geq 1\\0\leq r\leq N}}
        q^{\sum_i u_i(u_i-1)}\mathcal D_{a+1}(\mathbf u)
\]
be the shifted sum over upper chains of length at most $N$, and let
\[
        V_s\coloneqq
        \sum_{a\geq v_1\geq\cdots\geq v_s\geq 1}
        q^{\sum_i v_i(v_i+1)}\mathcal D_a(\mathbf v)
\]
be the contribution from lower chains of exact length $s$, with $V_0=1$. The
restricted contribution, after division by $(d+1)p^{da(a+1)}$, is at least
\[
        \qgauss{2a+1}{a}\frac{1}{d+1}
        \sum_{s=0}^{d-1} V_s U_{d-1-s}.
\]
Lemma~\ref{lem:shifted-chain} gives $U_N/(N+1)\to 1/(q;q)_a$ and
$U_N\leq (N+1)/(q;q)_a$. Also, by \eqref{eq:durfee-chain-shifted-square-sum},
\[
        \sum_{s\geq 0}V_s
        =
        \prod_{r=2}^{a+1}\frac{1}{1-q^r}
        =
        \frac{1-q}{(q;q)_{a+1}}.
\]
We record the limiting passage, since the factor $d+1$ is part of the leading
term. Put $L=(q;q)_a^{-1}$. For a fixed integer $M\geq 0$ and $d>M+1$,
\[
        \frac{1}{d+1}\sum_{s=0}^{d-1}V_sU_{d-1-s}
        \geq
        \sum_{s=0}^{M}V_s\frac{U_{d-1-s}}{d+1}.
\]
For each fixed $s\leq M$,
\[
        \frac{U_{d-1-s}}{d+1}
        =
        \frac{d-s}{d+1}\cdot \frac{U_{d-1-s}}{d-s}
        \longrightarrow L,
\]
since $d-s=(d-1-s)+1$. It follows that
\[
        \liminf_{d\to\infty}
        \frac{1}{d+1}\sum_{s=0}^{d-1}V_sU_{d-1-s}
        \geq
        L\sum_{s=0}^{M}V_s.
\]
Letting $M\to\infty$ and using the convergence of $\sum_sV_s$, we obtain the
lower limit
\[
        \qgauss{2a+1}{a}\frac{1}{(q;q)_a}\frac{1-q}{(q;q)_{a+1}}
        =\xi_{2a+1}(q).
\]
The upper bound gives the matching limsup. If $a=0$, sharpness is the exact
identity $s(C_{p^d})=d+1$.
\end{proofof}

\section{Block decoupling}
\label{sec:block-decoupling}

We now pass from a single homocyclic block to an arbitrary partition. The
columns of $\lambda$ of height $h$ form a block of width $d_h(\lambda)$; in
isolation, this block is the homocyclic group $C_{p^{d_h(\lambda)}}^{\,h}$.
The Ferrers diagram is not a product of its height blocks, since the entries of
$\mu'$ must remain nonincreasing across adjacent blocks. The next proposition
shows that, after normalisation by $p^{\eta(\lambda)}$, removing these boundary
conditions can only enlarge Delsarte's sum.

\begin{proposition}\label{prop:block-decoupling}
For every partition $\lambda$,
\[
        p^{-\eta(\lambda)}s(A_\lambda)
        \leq
        \prod_{h\geq 1}
        p^{-d_h(\lambda)\lfloor h^2/4\rfloor}
        s(C_{p^{d_h(\lambda)}}^{\,h}),
\]
where a factor with $d_h(\lambda)=0$ is interpreted as $1$.
\end{proposition}

\begin{proofof}
By \eqref{eq:normalised-delsarte-summand},
\[
        p^{-\eta(\lambda)}\alpha_\lambda(\mu;p)
        =
        q^{\sum_{i\geq 1}(\lfloor(\lambda'_i)^2/4\rfloor
        -\mu'_i(\lambda'_i-\mu'_i))}
        \prod_{i\geq 1}
        \qgauss{\lambda'_i-\mu'_{i+1}}{\mu'_i-\mu'_{i+1}}.
\]
Group the indices for which $\lambda'_i=h$. Such a block has length
$d=d_h(\lambda)$. Inside a fixed block write
$x_1\geq x_2\geq\cdots\geq x_d$ for the values of $\mu'$ in that block, and
let $y\leq x_d$ be the next value of $\mu'$ after the block; if no later block
is present, put $y=0$. The $q$-power part belonging to the block is already
internal. Writing $x_{d+1}=y$, the only boundary term in the Gaussian product is
the last factor of
\[
        \prod_{j=1}^{d}
        \qgauss{h-x_{j+1}}{x_j-x_{j+1}}.
\]
We compare this with the product obtained by replacing the lower boundary value
$y$ by $0$. Only the last factor changes, and the ratio of the new last factor
to the old one is
\[
        \frac{\qgauss{h}{x_d}}{\qgauss{h-y}{x_d-y}}
        =
        \frac{(q;q)_h(q;q)_{x_d-y}}{(q;q)_{h-y}(q;q)_{x_d}}
        =
        \prod_{k=1}^{y}
        \frac{1-q^{h-y+k}}{1-q^{x_d-y+k}}.
\]
For $y=0$ this is the ratio of two empty products. For $y>0$, the inequalities
$h\geq x_d$ and $0<q<1$ show that each numerator factor is at least the
corresponding denominator factor, whence the ratio is at least $1$. A
normalised summand is therefore not decreased when the height block is cut off from the
blocks below it by replacing its lower boundary value by $0$.

We may now perform this replacement for every height block. The monotonicity
conditions between different blocks may then be dropped; since every normalised
summand is a product of a power of $q$ and $q$-Gaussian coefficients, and so is
nonnegative, dropping these conditions only enlarges the sum. A block of height $h$ and width
$d=d_h(\lambda)$ contributes, with $x_{d+1}=0$,
\[
 \sum_{h\geq x_1\geq\cdots\geq x_d\geq 0}
 q^{\sum_{j=1}^{d}(\lfloor h^2/4\rfloor-x_j(h-x_j))}
 \prod_{j=1}^{d}
 \qgauss{h-x_{j+1}}{x_j-x_{j+1}}.
\]
By the homocyclic identity used in the proof of
Proposition~\ref{prop:homocyclic}, this is
$p^{-d\lfloor h^2/4\rfloor}s(C_{p^d}^{\,h})$. Multiplying over the independent
height blocks gives the claimed inequality.
\end{proofof}

We are now in a position to prove Theorem~\ref{thm:main}.

\begin{proofof}
The lower bound and the asymptotic for fixed $\lambda$ were proved in
Section~\ref{sec:delsarte}. For the upper bound, combine
Proposition~\ref{prop:block-decoupling} with
Proposition~\ref{prop:homocyclic}. A block of height one contributes exactly
$d_1(\lambda)+1$, the number of subgroups of the cyclic group
$C_{p^{d_1(\lambda)}}$. A nonzero block of even height $h=2a\geq 2$ contributes
at most $\xi_{2a}(q)$, and one of odd height $h=2a+1\geq 3$ at most
$(d_h(\lambda)+1)\xi_{2a+1}(q)$.

Multiplying the block estimates over all heights, the factors
$d_h(\lambda)+1$ from the blocks of odd height, including the block of height one,
give
\[
        \prod_{\substack{h\geq 1\\ h\textnormal{ odd}}}(d_h(\lambda)+1)
        =\theta(\lambda).
\]
The remaining factors are precisely the constants $\xi_h(q)$ for those heights
$h\geq 2$ with $d_h(\lambda)>0$. Restoring the factor $p^{\eta(\lambda)}$ in
Proposition~\ref{prop:block-decoupling} gives the asserted upper bound.
\end{proofof}

\begin{corollary}\label{cor:uniform}
For every partition $\lambda$,
\[
        s(A_\lambda)
        \leq
        K_p^{3r(\lambda)}\theta(\lambda)p^{\eta(\lambda)},
        \quad \textnormal{where } K_p\coloneqq(q;q)_\infty^{-1}.
\]
In particular, $s(A_\lambda)\leq 3^{r(\lambda)}\theta(\lambda)p^{\eta(\lambda)}$
for $p\geq 5$, and
$s(A_\lambda)\leq 6^{r(\lambda)}\theta(\lambda)3^{\eta(\lambda)}$ for $p=3$.
\end{corollary}

\begin{proofof}
It is enough to bound each block constant by $K_p^3$. Since the finite products
$(q;q)_a$ decrease to $(q;q)_\infty$, we have
$(q;q)_a\geq (q;q)_\infty$ for every finite $a$. If $h=2a$, then
\[
        (q;q)_{2a}
        =
        (q;q)_a\prod_{j=a+1}^{2a}(1-q^j)
        \leq (q;q)_a,
\]
and hence
\[
        \xi_{2a}(q)=
        \frac{(q;q)_{2a}}{(q;q)_a^4}
        \leq
        \frac{1}{(q;q)_a^3}
        \leq
        K_p^3.
\]
If $h=2a+1\geq 3$, then
\[
\begin{aligned}
        \xi_{2a+1}(q)
        &=(1-q)\frac{(q;q)_{2a+1}}{(q;q)_a^2(q;q)_{a+1}^2}  \\
        &\leq
        \frac{1}{(q;q)_a^2(q;q)_{a+1}}
        \leq K_p^3,
\end{aligned}
\]
because $(1-q)(q;q)_{2a+1}/(q;q)_{a+1}$ is a product of factors between $0$
and $1$. Inserting these estimates into Theorem~\ref{thm:main} proves the first
assertion.

It remains to justify the numerical constants. The product $K_p$ decreases with
$p$, since each factor $(1-p^{-j})^{-1}$ decreases with $p$; it therefore
suffices to bound $K_5^3$ for all $p\geq 5$, and $K_3^3$ for $p=3$. We use
$-\log(1-x)\leq x/(1-x)$ for $0<x<1$. Hence, for every $N\geq 1$,
\[
\begin{aligned}
        K_p^3
        &=
        \prod_{j=1}^{N}(1-p^{-j})^{-3}
        \exp\left(3\sum_{j>N}-\log(1-p^{-j})\right)       \\
        &\leq
        \prod_{j=1}^{N}(1-p^{-j})^{-3}
        \exp\left(3\sum_{j>N}\frac{p^{-j}}{1-p^{-j}}\right)       \\
        &\leq
        \prod_{j=1}^{N}(1-p^{-j})^{-3}
        \exp\left(
        \frac{3p^{-(N+1)}}{(1-p^{-1})(1-p^{-(N+1)})}
        \right).
\end{aligned}
\]
In the last line we used $1-p^{-j}\geq 1-p^{-(N+1)}$ for $j>N$, and then
summed the resulting geometric series. 
Taking $N=1$ gives
\[
        K_5^3
        \leq
        \left(5/4\right)^3
        \exp\left(5/32\right)
        <3,
\]
and
\[
        K_3^3
        \leq
        \left(3/2\right)^3
        \exp\left(9/16\right)
        <6.
\]
This proves the stated numerical bounds.
\end{proofof}

\begin{remark}
The case $r(\lambda)=1$ is particularly transparent. In this case there is a
unique height $h\geq 2$ for which $d_h(\lambda)>0$. Equivalently,
\[
        \lambda'=(h^a,1^b)
\]
for some $h\geq 2$, $a\geq 1$, and $b\geq 0$, or, in terms of $\lambda$,
\[
        \lambda=(a+b,\underbrace{a,\ldots,a}_{h-1\textnormal{ times}}).
\]
Theorem~\ref{thm:main} and Corollary~\ref{cor:uniform} give
\[
        \theta(\lambda)p^{\eta(\lambda)}
        \leq s(A_\lambda)
        \leq K_p^3\theta(\lambda)p^{\eta(\lambda)}.
\]
Since $K_3^3<6$, one has
\[
        \theta(\lambda)p^{\eta(\lambda)}
        \leq s(A_\lambda)
        <6\theta(\lambda)p^{\eta(\lambda)}
\]
for every odd prime $p$. For $p\geq 5$, the sharper estimate $K_5^3<3$
gives
\[
        \theta(\lambda)p^{\eta(\lambda)}
        \leq s(A_\lambda)
        <3\theta(\lambda)p^{\eta(\lambda)}.
\]
Thus, when only one nontrivial height block occurs, the leading term controls the whole abelian
subgroup count up to an absolute factor smaller than $6$, and up to a factor
smaller than $3$ for $p\geq 5$.
\end{remark}

\subsection{How large can the shape factors be?}

Theorem~\ref{thm:main} separates the main exponential term
$p^{\eta(\lambda)}$ from the remaining factors, 
namely the factor $\theta(\lambda)$ counting central patterns 
and the product of homocyclic block constants.
We record a coarse bound showing that these are small on the
exponential scale of $|\lambda|$.

Recall that $d_h(\lambda)$ is the number of columns of height exactly $h$. Thus
\[
        |\lambda|=\sum_{h\geq 1} h d_h(\lambda).
\]
If $r(\lambda)=r$, then there are $r$ distinct heights at least $2$ for which
$d_h(\lambda)>0$. The least possible contribution to $|\lambda|$ from such
heights is obtained by taking the heights $2,3,\ldots,r+1$, each with
multiplicity one. Hence
\[
        |\lambda|\geq 2+3+\cdots+(r+1)=\frac{r(r+3)}{2}.
\]
Solving $r(r+3)\leq 2|\lambda|$ gives
\[
        r(\lambda)
        \leq
        \left\lfloor \frac{\sqrt{8|\lambda|+9}-3}{2}\right\rfloor.
\]

We also need to control $\theta(\lambda)$. Define
\[
        \psi(u)\coloneqq
        \begin{cases}
        \log(1/u)-1+u, & 0<u\leq 1,\\
        0, & u\geq 1.
        \end{cases}
\]
For every $u>0$ and $x\geq 0$,
\[
        \log(x+1)\leq ux+\psi(u).
\]
Indeed, the maximum of $\log(x+1)-ux$ over $x\geq 0$ is $0$ if $u\geq 1$,
and is attained at $x=1/u-1$ with value $\log(1/u)-1+u$ if $0<u\leq 1$.
Taking $u=\tau h$, with $\tau>0$, and $x=d_h(\lambda)$, and summing over odd
$h$, gives
\[
        \log\theta(\lambda)
        \leq
        \tau\sum_{\substack{h\geq 1\\ h\textnormal{ odd}}}h d_h(\lambda)
        +
        \sum_{\substack{h\geq 1\\ h\textnormal{ odd}}}\psi(\tau h)
        \leq
        \tau|\lambda|+
        \sum_{\substack{h\geq 1\\ h\textnormal{ odd}}}\psi(\tau h).
\]

It remains to estimate the last sum. On $(0,1)$ we have $\psi''(u)=1/u^2>0$;
on $[1,\infty)$ the function is identically zero; and the one-sided
derivatives at $u=1$ are both $0$, so that $\psi$ is convex on $(0,\infty)$. Applying
Jensen's inequality to the intervals $[2k\tau,2(k+1)\tau]$, whose midpoints
are $(2k+1)\tau$, gives
\[
        2\tau\,\psi((2k+1)\tau)
        \leq
        \int_{2k\tau}^{2(k+1)\tau} \psi(u)\,du.
\]
For $k=0$ the endpoint $0$ is harmless; one may apply Jensen on $(0,2\tau)$
and then pass to the limit. Summing over $k\geq 0$ yields
\[
        \sum_{k\geq 0} \psi((2k+1)\tau)
        \leq
        \frac{1}{2\tau}\int_0^\infty \psi(u)\,du.
\]
Since $\psi(u)=0$ for $u\geq 1$ and
\[
        \int_0^\infty \psi(u)\,du
        =
        \int_0^1\bigl(\log(1/u)-1+u\bigr)\,du
        =
        1-1+\frac12
        =
        \frac12,
\]
we obtain
\[
        \sum_{\substack{h\geq 1\\ h\textnormal{ odd}}}\psi(\tau h)
        =
        \sum_{k\geq 0} \psi((2k+1)\tau)
        \leq
        \frac{1}{4\tau}.
\]
Therefore
\[
        \log\theta(\lambda)\leq \tau|\lambda|+\frac{1}{4\tau}.
\]
If $|\lambda|=0$, then $\theta(\lambda)=1$. If $|\lambda|\geq 1$, the right-hand side is minimised at $\tau=1/(2\sqrt{|\lambda|})$, and we get
\[
        \theta(\lambda)\leq \exp\sqrt{|\lambda|}.
\]

Combining this estimate with Corollary~\ref{cor:uniform} and the preceding
bound for $r(\lambda)$ gives
\[
        s(A_\lambda)
        \leq
        p^{\eta(\lambda)}
        \exp\left(
        \sqrt{|\lambda|}
        +
        3\log K_p
        \left\lfloor \frac{\sqrt{8|\lambda|+9}-3}{2}\right\rfloor
        \right).
\]
In particular, for fixed $p$,
\[
        s(A_\lambda)
        \leq
        p^{\eta(\lambda)}
        \exp\left((1+3\sqrt{2}\log K_p)\sqrt{|\lambda|}\right),
\]
since
$\lfloor(\sqrt{8|\lambda|+9}-3)/2\rfloor\leq \sqrt{2|\lambda|}$ for
$|\lambda|\geq 0$.

\section{Examples and sharpness}
\label{sec:examples-sharpness}

We record a few examples and explain why the uniform block constant $K_p^3$
cannot be improved. If $r(\lambda)=0$, then $A_\lambda$ is cyclic, say
$A_\lambda=C_{p^e}$ with $e\geq 0$. Here $\eta(\lambda)=0$ and
$\theta(\lambda)=e+1$, and Theorem~\ref{thm:main} gives the exact identity
$s(C_{p^e})=e+1$.

The condition $\theta(\lambda)=1$ is also easy to read from the diagram.
Since every factor in the definition of $\theta(\lambda)$ is a positive
integer, it holds if and only if $d_h(\lambda)=0$ for every odd $h$; that is,
exactly when every column of the Ferrers diagram has even height, so that the
rows of the diagram occur in equal pairs,
\[
        \lambda_1=\lambda_2,
        \quad
        \lambda_3=\lambda_4,
        \quad
        \lambda_5=\lambda_6,
        \quad\ldots ,
\]
and every part of $\lambda$ occurs with even multiplicity. 

The simplest noncyclic example is $C_{p^a}\times C_{p^a}$, with $a\geq 1$, 
for which $\lambda=(a,a)$, $\eta(\lambda)=a$, and $\theta(\lambda)=1$.

For homocyclic groups, Proposition~\ref{prop:homocyclic} gives the exact
limiting constants
\[
        \lim_{d\to\infty}p^{-da^2}s(C_{p^d}^{\,2a})
        =
        \frac{(q;q)_{2a}}{(q;q)_a^4}
        \qquad (a\geq 1)
\]
and
\[
        \lim_{d\to\infty}
        \frac{s(C_{p^d}^{\,2a+1})}{(d+1)p^{da(a+1)}}
        =
        (1-q)\frac{(q;q)_{2a+1}}{(q;q)_a^2(q;q)_{a+1}^2}
        \qquad (a\geq 0).
\]
These constants explain the size of the uniform block constant in
Corollary~\ref{cor:uniform}.
For the even-height block constants,
\[
        \frac{(q;q)_{2a}}{(q;q)_a^4}
        =
        \frac{\prod_{j=a+1}^{2a}(1-q^j)}{(q;q)_a^3}.
\]
The product in the numerator tends to $1$, since
\[
        \frac{(q;q)_\infty}{(q;q)_a}
        =
        \prod_{j>a}(1-q^j)
        \leq
        \prod_{j=a+1}^{2a}(1-q^j)
        \leq 1,
\]
and $(q;q)_a\to(q;q)_\infty$ as $a\to\infty$. Hence
\[
        \frac{(q;q)_{2a}}{(q;q)_a^4}
        \longrightarrow
        (q;q)_\infty^{-3}=K_p^3.
\]
This also shows that the uniform constant $K_p^3$ is best possible for a
single homocyclic block. 
If $C<K_p^3$, then $\xi_{2a}(q)>C$ for some $a$. 
Since the even homocyclic estimate in Proposition~\ref{prop:homocyclic} 
is sharp as $d\to\infty$, we have $s(C_{p^d}^{\,2a})>C p^{da^2}$ for all sufficiently large $d$. 
Thus no constant smaller than $K_p^3$ can replace it in a uniform estimate for a
single nonzero homocyclic block of arbitrary height and width. 
The odd-height constant, by contrast, tends to the strictly smaller limit
$(1-q)K_p^3$, since $0<q<1$ and $(q;q)_n\to(q;q)_\infty$ as $n\to\infty$; 
the blocks of even height are therefore the ones that pin down 
the value $K_p^3$ in Corollary~\ref{cor:uniform}.

For the numerical table below, $p=5$, and $B_\lambda$ denotes the refined upper
bound in Theorem~\ref{thm:main}:
\[
        B_\lambda\coloneqq
        \theta(\lambda)5^{\eta(\lambda)}
        \prod_{\substack{h\geq 2\\ d_h(\lambda)>0}}\xi_h(1/5).
\]
The exact values $s(A_\lambda)$ in the table were obtained by summing
Proposition~\ref{prop:delsarte} over all subpartitions $\mu\subseteq\lambda$;
the logarithms and ratios are rounded to three decimal places.

\begin{table}[htbp]
\centering
\begin{tabular}{@{}lrrrrrr@{}}
\toprule
Type $\lambda$ & $\eta(\lambda)$ & $\theta(\lambda)$ & $r(\lambda)$
& $\log_5 s(A_\lambda)$ & $\log_5 B_\lambda$ & $B_\lambda/s(A_\lambda)$ \\
\midrule
$(30)$ & $0$ & $31$ & $0$ & $2.134$ & $2.134$ & $1.000$ \\
$(24,6)$ & $6$ & $19$ & $1$ & $7.984$ & $8.220$ & $1.462$ \\
$(20,5^2)$ & $10$ & $96$ & $1$ & $12.962$ & $13.134$ & $1.318$ \\
$(18,6^2)$ & $12$ & $91$ & $1$ & $14.935$ & $15.100$ & $1.306$ \\
$(16,8,6)$ & $14$ & $63$ & $2$ & $16.745$ & $17.263$ & $2.301$ \\
$(15^2)$ & $15$ & $1$ & $1$ & $15.391$ & $15.391$ & $1.000$ \\
$(15,10,5)$ & $15$ & $36$ & $2$ & $17.415$ & $17.915$ & $2.236$ \\
$(12,10,8)$ & $18$ & $27$ & $2$ & $20.279$ & $20.736$ & $2.086$ \\
$(10^3)$ & $20$ & $11$ & $1$ & $21.764$ & $21.788$ & $1.039$ \\
$(12,8,6,4)$ & $22$ & $15$ & $3$ & $24.110$ & $24.857$ & $3.325$ \\
$(9,8,7,6)$ & $27$ & $4$ & $3$ & $28.381$ & $29.036$ & $2.866$ \\
$(8,7,6,5,4)$ & $31$ & $20$ & $4$ & $33.359$ & $34.393$ & $5.284$ \\
$(6^5)$ & $36$ & $7$ & $1$ & $37.522$ & $37.566$ & $1.074$ \\
$(6,5^4,4)$ & $42$ & $4$ & $2$ & $43.351$ & $43.724$ & $1.823$ \\
$(5^6)$ & $45$ & $1$ & $1$ & $45.505$ & $45.506$ & $1.001$ \\
$(4^7,2)$ & $56$ & $3$ & $2$ & $57.176$ & $57.561$ & $1.860$ \\
$(4^5,2^5)$ & $62$ & $3$ & $2$ & $63.397$ & $63.550$ & $1.279$ \\
$(3^{10})$ & $75$ & $1$ & $1$ & $75.502$ & $75.511$ & $1.014$ \\
$(2^8,1^{14})$ & $137$ & $1$ & $2$ & $137.761$ & $138.020$ & $1.519$ \\
$(1^{30})$ & $225$ & $1$ & $1$ & $225.381$ & $225.511$ & $1.233$ \\
\bottomrule
\end{tabular}
\vspace{0.3em}
\caption{Comparison, for selected types $\lambda$ with $p=5$ and
$|\lambda|=30$, between the exact value $s(A_\lambda)$ and the refined bound
$B_\lambda$.}
\label{tab:comparison-30}
\end{table}

\FloatBarrier
\section*{Applications}

Sections~\ref{sec:overview}--\ref{sec:examples-sharpness} prove the abelian estimate and discuss its sharpness. The remaining
sections use the bound in three directions. We first transfer it to groups of
nilpotency class less than $p$, and then return to abelian groups in two
asymptotic settings.

\section{Groups of nilpotency class less than \texorpdfstring{$p$}{p}}
\label{sec:class-less-than-p}

We begin with the application to groups of nilpotency class less than $p$.
By a theorem of Groves, such a group admits an abelian group law on the same
underlying set, in such a way that every subgroup of the original group
remains a subgroup of the abelian one; subgroup counting is thereby reduced
to the abelian case.

The comparison group is determined by the following logarithmic
$\Omega$-profile. If $G$ is a finite
$p$-group, let $\Omega_i(G)$ be the subgroup generated by the elements of
order dividing $p^i$, with $\Omega_0(G)=1$. If $\exp(G)=p^e$, put, for
$0\leq i\leq e$,
\[
        \omega_i(G)=\log_p\card{\Omega_i(G)}.
\]
When $G$ is fixed, write $\omega_i$ for $\omega_i(G)$. The sequence
\[
        0=\omega_0\leq \omega_1\leq\cdots\leq \omega_e=\log_p\card G
\]
is the logarithmic $\Omega$-profile of $G$, abbreviated by the tuple
$\omega(G)=(\omega_0,\omega_1,\ldots,\omega_e)$.

\begin{theorem}\label{thm:groves-reduction}
Let $p$ be a prime, and let $G$ be a finite $p$-group of nilpotency class
less than $p$. Write $\exp(G)=p^e$, set $\omega_0=0$, and put
$\omega_i=\log_p\card{\Omega_i(G)}$ for $1\leq i\leq e$. Then
\[
        s(G)\leq s(A_\lambda),
\]
where $A_\lambda$ is the abelian $p$-group whose type satisfies, for
$1\leq i\leq e$,
\[
        \lambda'_i=\omega_i-\omega_{i-1}.
\]
Equivalently,
\[
        A_\lambda\cong \prod_{k=1}^e C_{p^k}^{w_k},
\]
where
\[
        w_k=
        \begin{cases}
        2\omega_k-\omega_{k+1}-\omega_{k-1}, & k<e,\\
        \omega_e-\omega_{e-1}, & k=e.
        \end{cases}
\]
\end{theorem}

\begin{proofof}
By Groves~\cite[Th.~2]{Groves}, $G$ is verbally abelian; there is
a group word $v(x,y)=xyw(x,y)$, with $w$ a commutator word, such that
\[
        x*y=v(x,y)
\]
makes the underlying set of $G$ into an abelian group. Let $A=(G,*)$. We
claim that every subgroup of $G$ is a subgroup of $A$. Indeed, if
$H\leq G$ and $x,y\in H$, then $w(x,y)\in H$, since $w$ is a word in
commutators of $x$ and $y$. Hence $x*y=xyw(x,y)\in H$. The identity element is
the same for the two operations. Moreover, the operations agree on every cyclic
subgroup, because commutator words vanish on commuting elements. Thus inverses
and element orders are the same in $G$ and in $A$. It follows that every
subgroup of $G$ is a subgroup of $A$, whence $s(G)\leq s(A)$.

Next, a $p$-group of class less than $p$ is regular, and in a regular
$p$-group the subgroup $\Omega_i(G)$ is exactly the set of elements of order
dividing $p^i$; see Huppert~\cite[Kap.~III, \S 10, Satz~10.2 and
Hauptsatz~10.5]{Huppert}. For the abelian group $A$, the same description holds
directly. Since $G$ and $A$ have the same element orders, for $0\leq i\leq e$ we have
\[
        \card{\Omega_i(A)}=\card{\Omega_i(G)}=p^{\omega_i}.
\]

Write
\[
        A\cong \prod_{k=1}^e C_{p^k}^{w_k}.
\]
If $\lambda$ is the type of $A$, then
\[
        \lambda'_i=\sum_{k\geq i}w_k.
\]
On the other hand,
\[
        \omega_i=\log_p\card{\Omega_i(A)}
        =
        \sum_{k=1}^e w_k\min\{i,k\}.
\]
Subtracting the same formula with $i-1$ in place of $i$ gives
\[
        \omega_i-\omega_{i-1}
        =
        \sum_{k\geq i}w_k
        =
        \lambda'_i.
\]
The increments $\omega_i-\omega_{i-1}$ thus form the conjugate type $\lambda'$
of the abelian group $A$, and in particular constitute a partition. Solving
for the multiplicities $w_k$ gives
\[
        w_e=\omega_e-\omega_{e-1}
\]
and, for $k<e$,
\[
        w_k
        =
        (\omega_k-\omega_{k-1})-(\omega_{k+1}-\omega_k)
        =
        2\omega_k-\omega_{k+1}-\omega_{k-1},
\]
as claimed.
\end{proofof}

The same conclusion can also be reached from Lazard's correspondence. Indeed,
$G$ corresponds to a Lie ring $L$ on the same underlying set, and the subgroups
of $G$ correspond to Lie subrings of $L$. Every Lie subring is an additive
subgroup of $(L,+)$, so $s(G)\leq s(L,+)$. Groves's theorem gives the route
used here more directly; see Lazard~\cite[Chap.~II, Th.~(4.6)]{Lazard}.

\begin{corollary}\label{cor:omega-profile-bound}
Let $G$ be a finite $p$-group of nilpotency class less than $p$, and let
$\lambda$ be the partition of Theorem~\ref{thm:groves-reduction}, so that
$\lambda'_i=\omega_i-\omega_{i-1}$ for $1\leq i\leq e$. Then
\[
        s(G)
        \leq
        \theta(\lambda)p^{\eta(\lambda)}
        \prod_{\substack{h\geq 2\\ d_h(\lambda)>0}}\xi_h(p^{-1}).
\]
In particular,
\[
        s(G)
        \leq
        K_p^{3r(\lambda)}\theta(\lambda)p^{\eta(\lambda)}.
\]
\end{corollary}

\begin{proofof}
By Theorem~\ref{thm:groves-reduction}, $s(G)\leq s(A_\lambda)$; the two
displayed bounds are those of Theorem~\ref{thm:main} and
Corollary~\ref{cor:uniform} for $A_\lambda$.
\end{proofof}

\begin{remark}
The estimate in Corollary~\ref{cor:omega-profile-bound} depends only on the
$\Omega$-profile of $G$. The commutator structure enters only through the class
assumption, which supplies the comparison with an abelian group. Thus two groups
with the same $\Omega$-profile receive the same bound, even when their subgroup
lattices are different.

There are two boundary cases in which the comparison gives no new information.
If $p=2$, the class condition forces $G$ to be abelian. If $\exp(G)=p$, the
$\Omega$-profile records only $\card G$, and the associated abelian group is
elementary abelian.
\end{remark}

\subsection{Extraspecial groups of exponent \texorpdfstring{$p^2$}{p2}}

We next test the sharpness of the abelian comparison on extraspecial groups. Throughout this
subsection, $p$ is odd and $G_n$ is an extraspecial $p$-group of order
$p^{2n+1}$ and exponent $p^2$. Put
\[
        Z=Z(G_n)=G_n'=\Phi(G_n).
\]
Then $Z\cong C_p$, and $V=G_n/Z$ is a $2n$-dimensional vector space over
$\mathbb F_p$. The equality $Z=G_n'=\Phi(G_n)$ is part of the defining
structure of extraspecial groups; see Huppert~\cite[Kap.~III, \S 13,
Def.~13.1]{Huppert}. We shall use the alternating form afforded by the
commutator map; compare Bouc--Mazza~\cite[Sec.~2]{BoucMazza}. Our notation
for the symplectic space is standard; see Taylor~\cite[Chap.~8]{Taylor}.

Choose a generator $z$ of $Z$.
The commutator gives a nondegenerate alternating bilinear form
\[
        \beta:V\times V\longrightarrow \mathbb F_p,
        \qquad
        [x,y]=z^{\beta(xZ,yZ)}.
\]
Since $G_n/Z$ has exponent $p$, we have $x^p\in Z$ for every $x\in G_n$.
Moreover, if $i\in\mathbb Z$, then $(xz^i)^p=x^p$, since $z$ has order $p$
and is central. Hence there is a well-defined map
\[
        f:V\longrightarrow \mathbb F_p,
        \qquad
        x^p=z^{f(xZ)}.
\]
The collection formula in class two gives
\[
        (xy)^p=x^py^p[y,x]^{\binom p2}=x^py^p,
\]
because $[y,x]$ has order $p$ and $p\mid \binom p2$. Thus $f$ is additive, and
therefore $\mathbb F_p$-linear. Since $G_n$ has exponent $p^2$, this linear
form is nonzero, so $\ker f$ is a hyperplane of $V$. Therefore
\[
        \Omega_1(G_n)=\{x\in G_n:x^p=1\}
        =\pi^{-1}(\ker f),
\]
where $\pi:G_n\to V$ is the quotient map, and hence
\[
        \card{\Omega_1(G_n)}=p\cdot p^{2n-1}=p^{2n}.
\]
As $\Omega_2(G_n)=G_n$, the associated abelian comparison group is
\[
        A_n\cong C_{p^2}\times C_p^{\,2n-1}.
\]
Indeed, here $\omega_1=2n$, $\omega_2=2n+1$, and Theorem~\ref{thm:groves-reduction}
gives $w_2=1$ and $w_1=2n-1$.

We now count the subgroups of $G_n$. Let
\[
        T_n=\sum_{j=0}^{2n}\gauss{2n}{j},
\]
the number of subspaces of $V$. The subgroups containing $Z$ are precisely the
preimages of subspaces of $V$, and therefore contribute $T_n$ subgroups.

It remains to count the subgroups $H\leq G_n$ with $H\cap Z=1$. Such a
subgroup is elementary abelian, since $h^p\in H\cap Z$ and
$[h,k]\in H\cap Z$ for $h,k\in H$. If $U=HZ/Z\leq V$, then $H\to U$ is
an isomorphism, and $U$ is contained in $\ker f$ and is totally isotropic
for $\beta$. Conversely, if $U\leq\ker f$ is a totally isotropic subspace of
dimension $r$, then $\pi^{-1}(U)$ is elementary abelian of order $p^{r+1}$
and contains $Z$ as a one-dimensional subspace. The complements to $Z$ in
this vector space are the graphs of the linear maps $U\to Z$, so there are
$p^r$ of them.

Let $I_m(t)$ be the number of $t$-dimensional totally isotropic subspaces of
a nondegenerate symplectic space of dimension $2m$ over $\mathbb F_p$, with
$I_m(t)=0$ if $t<0$ or $t>m$; in particular, $I_m(0)=1$. We obtain $I_m(t)$ by the usual count of ordered bases. After $i$ independent isotropic vectors have been chosen, 
their span has dimension $i$ and its perpendicular space has dimension $2m-i$, 
so the next vector has $p^{2m-i}-p^i$ choices. Dividing by the number
$\prod_{i=0}^{t-1}(p^t-p^i)$ of ordered bases of a fixed subspace of
dimension $t$ gives
\[
        I_m(t)
        =
        \prod_{i=0}^{t-1}\frac{p^{2m-i}-p^i}{p^t-p^i}
        =
        \gauss{m}{t}\prod_{i=0}^{t-1}(p^{m-i}+1).
\]
Since $\beta$ is nondegenerate and $f$ is nonzero, there is a nonzero vector
$v_0\in V$ such that $f(v)=\beta(v,v_0)$ for every $v\in V$. Thus
$\ker f=v_0^\perp$, where $v_0^\perp=\{v\in V:\beta(v,v_0)=0\}$. Also
$v_0\in v_0^\perp$, since $\beta$ is alternating. If
$w\in\operatorname{rad}(v_0^\perp)$, then $w$ is orthogonal to all of
$v_0^\perp$, so $w\in(v_0^\perp)^\perp=\langle v_0\rangle$. Hence
\[
        \operatorname{rad}(\ker f)=\langle v_0\rangle.
\]
The quotient $\ker f/\langle v_0\rangle$ is therefore a nondegenerate
symplectic space of dimension $2(n-1)$.

We now count the $r$-dimensional totally isotropic subspaces of $\ker f$.
Those containing $\langle v_0\rangle$ correspond to $(r-1)$-dimensional totally
isotropic subspaces of $\ker f/\langle v_0\rangle$, and contribute
$I_{n-1}(r-1)$. If $U$ meets $\langle v_0\rangle$ trivially, then its image
$\overline U$ in the quotient has dimension $r$ and is totally isotropic. For
such a subspace $\overline U$, the inverse image in $\ker f$ is an extension
of $\overline U$ by $\langle v_0\rangle$, and the complements to
$\langle v_0\rangle$ in it are the graphs of the linear maps
$\overline U\to\langle v_0\rangle$. There are $p^r$ such maps, whence the number
of $r$-dimensional totally isotropic subspaces of $\ker f$ is
\[
        I_{n-1}(r-1)+p^rI_{n-1}(r).
\]

Thus we have the exact formula
\begin{equation}\label{eq:extraspecial-count}
        s(G_n)
        =
        \sum_{j=0}^{2n}\gauss{2n}{j}
        +
        \sum_{r=0}^{n}
        p^r\left(I_{n-1}(r-1)+p^rI_{n-1}(r)\right).
\end{equation}
For the corresponding poset-theoretic description of subgroups avoiding the
Frattini subgroup, see Bouc--Mazza~\cite[Secs.~2 and~4]{BoucMazza}.

The abelian comparison group $A_n=C_{p^2}\times C_p^{\,2n-1}$ has subgroup
count
\begin{equation}\label{eq:abelian-extraspecial-comparison}
        s(A_n)
        =
        \sum_{r=0}^{2n}\gauss{2n}{r}
        +
        \sum_{r=1}^{2n}p^{2n-r}\gauss{2n-1}{r-1}.
\end{equation}
The first sum counts the elementary abelian subgroups, while the second counts
those of type $C_{p^2}\times C_p^{r-1}$.

The inequality $s(G_n)\leq s(A_n)$ is the abelian comparison obtained above.
For the reverse estimate up to a factor $2$, use the Gaussian Pascal identity
\[
        \gauss{2n}{r}
        =
        \gauss{2n-1}{r}
        +
        p^{2n-r}\gauss{2n-1}{r-1}.
\]
It shows that the second sum in \eqref{eq:abelian-extraspecial-comparison} is
at most $T_n$. Since \eqref{eq:extraspecial-count} gives $T_n\leq s(G_n)$,
we obtain
\[
        s(G_n)\leq s(A_n)\leq 2T_n\leq 2s(G_n).
\]
Thus for extraspecial groups of exponent $p^2$, the abelian comparison loses
at most a factor $2$.

For fixed $n\geq 2$ and $p\to\infty$, the factor is best possible. Indeed, the
abelian group $A_n$ has two leading subgroup types,
\[
        C_p^n
        \qquad\textnormal{and}\qquad
        C_{p^2}\times C_p^{\,n-1},
\]
so
\[
        s(A_n)=2p^{n^2}+O_n(p^{n^2-1}).
\]
In $G_n$, the subgroups containing $Z$ contribute
$p^{n^2}+O_n(p^{n^2-1})$, and the second sum in
\eqref{eq:extraspecial-count} has degree $<n^2$. Indeed, the displayed product
for $I_m(t)$ gives
\[
        \deg_p I_m(t)=2mt-\frac{3t^2-t}{2}.
\]
Thus $p^{2r}I_{n-1}(r)$ has degree
$2nr-(3r^2-r)/2$, whose gap from $n^2$ is
$(n-r)^2+r(r-1)/2>0$ for $0\leq r\leq n-1$. The term
$p^rI_{n-1}(r-1)$ is smaller, and for $r=n$ it has degree
$n+n(n-1)/2<n^2$. Hence
\[
        s(G_n)=p^{n^2}+O_n(p^{n^2-1}),
\]
and therefore $s(A_n)/s(G_n)\to 2$.

\subsection{Cyclic extensions by power automorphisms}

The preceding example shows that the abelian comparison need not be exact. We
next record a family of nonabelian groups for which, by contrast, it is.

Let $p$ be odd. Let $B$ be an abelian $p$-group of exponent $p^a$, and
let $b\geq 1$, $1\leq \rho<a$, and $u$ prime to $p$. Put
\[
        \nu\coloneqq 1+up^\rho.
\]
Assume $\rho+b\geq a$. Then $\nu^{p^b}\equiv 1\pmod{p^a}$, so the map
$x\mapsto x^\nu$ defines an automorphism of $B$ of order dividing
$p^b$. Form the semidirect product
\[
        G=B\rtimes_\nu \langle y\rangle,
        \qquad
        \langle y\rangle\cong C_{p^b},
\]
where $yxy^{-1}=x^\nu$ for $x\in B$. Let
\[
        A=B\times C_{p^b}.
\]

We shall use the following elementary identity: for every integer $m\geq 0$
and every $t\geq 0$,
\begin{equation}\label{eq:geometric-sum-unit}
        1+\nu^m+\nu^{2m}+\cdots+\nu^{(p^t-1)m}
        =
        p^t\varepsilon
\end{equation}
for some integer $\varepsilon$ prime to $p$. If $m=0$, the left-hand side is
$p^t$. Suppose that $m\neq 0$. Write $v_p$ for the $p$-adic valuation,
normalized by $v_p(p)=1$. Since $\nu=1+up^\rho$, with $u$ prime to $p$, we have
\[
        v_p(\nu^m-1)=\rho+v_p(m).
\]
Indeed, write $m=p^e m_0$, with $p\nmid m_0$. The binomial expansion gives
$v_p((1+up^\rho)^{m_0}-1)=\rho$, because the first nonzero term is
$m_0up^\rho$ and all later terms have larger $p$-adic valuation. If
$\chi=1+cp^r$ with $p\nmid c$ and $r\geq 1$, then another binomial expansion
gives $v_p(\chi^p-1)=r+1$. Applying this step $e$ times gives
$v_p(\nu^m-1)=\rho+e$. Replacing $m$ by $mp^t$ gives
\[
        v_p((\nu^m)^{p^t}-1)=\rho+v_p(m)+t.
\]
Therefore
\[
        1+\nu^m+\cdots+\nu^{(p^t-1)m}
        =
        \frac{(\nu^m)^{p^t}-1}{\nu^m-1}
\]
has $p$-adic valuation $t$, which proves \eqref{eq:geometric-sum-unit}.

\begin{lemma}\label{lem:power-action-omega-profile}
For every $i\geq 0$,
\[
        \card{\Omega_i(G)}=\card{\Omega_i(A)}.
\]
Consequently $G$ and $A$ have the same $\Omega$-profile.
\end{lemma}

\begin{proofof}
Every element of $G$ has a unique expression $xy^m$, with $x\in B$ and
$0\leq m<p^b$, and
\[
        (xy^m)^{p^i}
        =
        x^{1+\nu^m+\nu^{2m}+\cdots+\nu^{(p^i-1)m}}
        y^{mp^i}.
\]
In view of \eqref{eq:geometric-sum-unit}, the first factor is trivial precisely when
$x^{p^i}=1$, and the second precisely when $y^{mp^i}=1$. Thus the elements
killed by the $p^i$-power map are exactly
\[
        \{xy^m:x^{p^i}=1,\ y^{mp^i}=1\}.
\]
The displayed set is closed under multiplication. Indeed, if
$x_1^{p^i}=x_2^{p^i}=1$, then $(x_2^{\nu^{m_1}})^{p^i}=1$, because
$\nu^{m_1}$ is prime to $p$, and the condition $y^{mp^i}=1$ is closed under
addition of the exponents of $y$. This set is therefore a subgroup, and hence
equal to $\Omega_i(G)$. Its cardinality is the same as that of
$\Omega_i(B\times C_{p^b})$.
\end{proofof}

\begin{proposition}\label{prop:power-action-subgroup-count}
With $G=B\rtimes_\nu C_{p^b}$ and $A=B\times C_{p^b}$ as above, we have
$s(G)=s(A)$. More precisely,
\[
        s(G)=s(A)=
        \sum_{U\leq B}\sum_{j=0}^{b}
        \card{\Omega_{b-j}(B/U)}.
\]
\end{proposition}

\begin{proofof}
Let $\pi:G\to G/B\cong C_{p^b}$ be the quotient map. For a subgroup
$H\leq G$, put $U=H\cap B$. The subgroup $\pi(H)$ is uniquely of the form
\[
        \langle y^{p^j}B\rangle
\]
for some $0\leq j\leq b$, with $j=b$ corresponding to the trivial image.

Since $y$ acts on $B$ by a power automorphism, every subgroup $U\leq B$ is
$y$-invariant, and hence normal in $G$. Fix $U\leq B$ and $0\leq j\leq b$, and
let $H\leq G$ satisfy $H\cap B=U$ and $\pi(H)=\langle y^{p^j}B\rangle$. Then
$H/U$ is cyclic, generated by the image of any element $xy^{p^j}\in H$ with
$x\in B$, so that $H=\langle U,xy^{p^j}\rangle$. Since $U$ is normal, this
subgroup equals $U\langle xy^{p^j}\rangle$, and its intersection with $B$ is
$U\langle (xy^{p^j})^{p^{b-j}}\rangle$. Hence
\[
        H\cap B=U
        \qquad\Longleftrightarrow\qquad
        (xy^{p^j})^{p^{b-j}}\in U.
\]
By \eqref{eq:geometric-sum-unit} the left-hand element equals
$x^{p^{b-j}\varepsilon}$ with $\varepsilon$ prime to $p$, so the condition is
equivalent to $x^{p^{b-j}}\in U$, that is,
\[
        xU\in \Omega_{b-j}(B/U).
\]
The assignment $H\mapsto xU$ is well defined: if $xy^{p^j}$ and $x'y^{p^j}$
both lie in $H$, then $(xy^{p^j})(x'y^{p^j})^{-1}=xx'^{-1}\in H\cap B=U$, so
the coset $xU$ does not depend on the chosen generator. It is a bijection from
the subgroups $H$ with $H\cap B=U$ and $\pi(H)=\langle y^{p^j}B\rangle$ onto
$\Omega_{b-j}(B/U)$: the inverse sends a coset $xU$ to
$\langle U,xy^{p^j}\rangle$, which is well defined because replacing $x$ by
another element of $xU$ alters $xy^{p^j}$ only by a factor in $U$, and which
has intersection $U$ with $B$ by the equivalence above. When $j=b$ one has
$\Omega_0(B/U)=1$, the unique subgroup being $U$ itself. Summing
$\card{\Omega_{b-j}(B/U)}$ over all $U\leq B$ and all $0\leq j\leq b$ gives the
displayed formula for $s(G)$.

The same argument applies to the direct product $A=B\times C_{p^b}$, where the
action is trivial. Therefore $s(G)=s(A)$.
\end{proofof}

\begin{corollary}\label{cor:power-action-sharpness}
Assume in addition that
\[
        \left\lceil\frac{a}{\rho}\right\rceil<p.
\]
Then $G=B\rtimes_\nu C_{p^b}$ has nilpotency class less than $p$, and
the abelian comparison in Theorem~\ref{thm:groves-reduction} is exact:
\[
        s(G)=s(A),
        \qquad
        A\cong B\times C_{p^b}.
\]
\end{corollary}

\begin{proofof}
For $x\in B$,
\[
        [y,x]=yxy^{-1}x^{-1}=x^{\nu-1}=x^{up^\rho}.
\]
Thus $\gamma_2(G)=B^{p^\rho}$. Suppose that $k\geq 2$ and
$\gamma_k(G)=B^{p^{(k-1)\rho}}$. Since $B$ is abelian, the only nontrivial
commutators with elements of $\gamma_k(G)$ are those with $y$. For $x\in B$,
\[
        [y,x^{p^{(k-1)\rho}}]=x^{up^{k\rho}}.
\]
Since $u$ is prime to $p$, this proves by induction that, for $k\geq 1$,
\[
        \gamma_{k+1}(G)=B^{p^{k\rho}}.
\]
Since $B$ has exponent $p^a$, the nilpotency class of $G$ is
$\lceil a/\rho\rceil$. The displayed hypothesis therefore gives class less
than $p$. In view of Lemma~\ref{lem:power-action-omega-profile}, the abelian group
associated with $G$ by Theorem~\ref{thm:groves-reduction} is $A$, and
Proposition~\ref{prop:power-action-subgroup-count} gives $s(G)=s(A)$.
\end{proofof}

\begin{remark}
Taking $B=C_{p^a}$ gives split metacyclic examples
\[
        C_{p^a}\rtimes_\nu C_{p^b},
        \qquad
        yxy^{-1}=x^\nu.
\]
The same construction applies to every abelian $p$-group $B$ of exponent
$p^a$. Thus it gives nonabelian groups of class less than $p$ for which the
abelian comparison is exact. For subgroup parameterisations of finite metacyclic
groups, see Yang~\cite[Th.~3.2]{Yang2020}.
\end{remark}

\subsection{Numerical examples}

We close this section by applying Theorem~\ref{thm:groves-reduction} to selected groups of order $5^9$. 
The tables below record the partitions obtained from their $\Omega$-profiles and the resulting ratios of subgroup counts.
We take $p=5$, so that the hypothesis of class less than $p$ allows the
nilpotency classes $2$, $3$, and $4$. For $p=3$ only class $2$ is possible
among nonabelian groups, and for $p=2$ the condition forces the group to be
abelian.
The groups are given by the power-commutator presentations in
Table~\ref{tab:comparison-5-9-presentations}. In the second column, an entry
$(a_1,\ldots,a_m)$ means that the displayed generators $x_1,\ldots,x_m$ have
orders $5^{a_1},\ldots,5^{a_m}$. In the third column we list precisely those commutators $[x_i,x_j]$ between
the displayed generators, with $i<j$, that are nontrivial. We used GAP~\cite{GAP4} to verify the generator orders, the group orders, and the nilpotency classes.

\begin{table}[htbp]
\centering
\setlength{\tabcolsep}{5pt}
\renewcommand{\arraystretch}{1.12}
\begin{tabular}{@{}cll@{}}
\toprule
Group & generator orders & nontrivial commutators\\
\midrule
$G_1$  & $(5,3,1)$       & $[x_1,x_2]=x_3$ \\
$G_2$  & $(4,4,1)$       & $[x_1,x_2]=x_3$ \\
$G_3$  & $(3,3,2,1)$     & $[x_1,x_2]=x_4$ \\
$G_4$  & $(3,3,2,1)$     & $[x_1,x_2]=x_3,\ [x_1,x_3]=x_4$ \\
$G_5$  & $(3,3,1,1,1)$   & $[x_1,x_2]=x_5$ \\
$G_6$  & $(4,2,1,1,1)$   & $[x_1,x_2]=x_5$ \\
$G_7$  & $(4,2,1,1,1)$   & $[x_1,x_2]=x_3,\ [x_1,x_3]=x_4,\ [x_1,x_4]=x_5$ \\
$G_8$  & $(5,1,1,1,1)$   & $[x_1,x_2]=x_5$ \\
$G_9$  & $(5,1,1,1,1)$   & $[x_1,x_2]=x_3,\ [x_1,x_3]=x_4,\ [x_1,x_4]=x_5$ \\
$G_{10}$ & $(3,3,3)$     & $[x_1,x_2]=x_3$ \\
$G_{11}$ & $(3,2,2,2)$   & $[x_1,x_2]=x_4$ \\
$G_{12}$ & $(2,2,1,1,3)$ & $[x_1,x_2]=x_3,\ [x_1,x_3]=x_4$ \\
$G_{13}$ & $(2,2,2,2,1)$ & $[x_1,x_2]=x_3,\ [x_1,x_3]=x_5$ \\
$G_{14}$ & $(2,2,2,2,1)$ & $[x_1,x_2]=x_3,\ [x_1,x_3]=x_4,\ [x_1,x_4]=x_5$ \\
\bottomrule
\end{tabular}
\vspace{0.3em}
\caption{Power-commutator presentations for the groups in
Table~\ref{tab:comparison-5-9}.}
\label{tab:comparison-5-9-presentations}
\end{table}

For each group $G_i$, the partition $\lambda$ in
Table~\ref{tab:comparison-5-9} is the type of the abelian comparison group
$A_\lambda$, obtained from the logarithmic $\Omega$-profile by
$\lambda'_j=\omega_j(G_i)-\omega_{j-1}(G_i)$. Here $c$ denotes the nilpotency
class.
The values of $s(G_i)$ were computed with GAP. 
For each displayed presentation, we constructed
the corresponding group, converted it to a polycyclic representation,
and enumerated its subgroups.

\begin{table}[htbp]
\centering
\setlength{\tabcolsep}{3.2pt}
\renewcommand{\arraystretch}{1.12}
\begin{tabular}{@{}cllrrrrrr@{}}
\toprule
$G$ & $\omega(G)$ & $\lambda$ & $c$ & $\eta(\lambda)$ & $\theta(\lambda)$ & $r(\lambda)$
& $s(G)$ & $s(A_\lambda)/s(G)$ \\
\midrule
$G_1$
& $(0,3,5,7,8,9)$ & $(5,3,1)$ & $2$ & $4$ & $6$ & $2$
& $5375$ & $1.004651$ \\

$G_2$
& $(0,3,5,7,9)$ & $(4,4,1)$ & $2$ & $5$ & $2$ & $2$
& $11625$ & $1.002151$ \\

$G_3$
& $(0,4,7,9)$ & $(3,3,2,1)$ & $2$ & $7$ & $2$ & $3$
& $357525$ & $1.052794$ \\

$G_4$
& $(0,4,7,9)$ & $(3,3,2,1)$ & $3$ & $7$ & $2$ & $3$
& $225000$ & $1.672889$ \\

$G_5$
& $(0,5,7,9)$ & $(3,3,1,1,1)$ & $2$ & $8$ & $2$ & $2$
& $1771400$ & $1.019758$ \\

$G_6$
& $(0,5,7,8,9)$ & $(4,2,1,1,1)$ & $2$ & $7$ & $6$ & $2$
& $771400$ & $1.045372$ \\

$G_7$
& $(0,5,7,8,9)$ & $(4,2,1,1,1)$ & $4$ & $7$ & $6$ & $2$
& $573125$ & $1.407023$ \\

$G_8$
& $(0,5,6,7,8,9)$ & $(5,1,1,1,1)$ & $2$ & $6$ & $10$ & $1$
& $171400$ & $1.204201$ \\

$G_9$
& $(0,5,6,7,8,9)$ & $(5,1,1,1,1)$ & $4$ & $6$ & $10$ & $1$
& $166125$ & $1.242438$ \\

$G_{10}$
& $(0,3,6,9)$ & $(3,3,3)$ & $2$ & $6$ & $4$ & $1$
& $32250$ & $2.803101$ \\

$G_{11}$
& $(0,4,8,9)$ & $(3,2,2,2)$ & $2$ & $8$ & $2$ & $1$
& $404400$ & $3.867211$ \\

$G_{12}$
& $(0,5,8,9)$ & $(3,2,2,1,1)$ & $3$ & $8$ & $8$ & $2$
& $884400$ & $5.151967$ \\

$G_{13}$
& $(0,5,9)$ & $(2,2,2,2,1)$ & $3$ & $10$ & $2$ & $2$
& $1740025$ & $21.745463$ \\

$G_{14}$
& $(0,5,9)$ & $(2,2,2,2,1)$ & $4$ & $10$ & $2$ & $2$
& $663750$ & $57.005876$ \\
\bottomrule
\end{tabular}
\vspace{0.3em}
\caption{Subgroup counts for the groups $G_i$ and comparison ratios.}
\label{tab:comparison-5-9}
\end{table}

Two features of Table~\ref{tab:comparison-5-9} are worth a remark. Among the
groups that share a comparison type, namely the pairs $(G_3,G_4)$,
$(G_6,G_7)$, $(G_8,G_9)$ and $(G_{13},G_{14})$, the ratio
$s(A_\lambda)/s(G)$ rises with the nilpotency class. This is consistent with
the idea that higher nilpotency class leaves more room for the subgroup lattice
to differ from that of the abelian comparison group. Across the table the
comparison ranges widely, from sharpness within one per cent for $G_1$ and
$G_2$ to an overshoot by a factor above $57$ for $G_{14}$. It would be
interesting to identify the structural features of $G$ that govern this ratio.

\FloatBarrier

\section{Fixed-shape diagonal sums}
\label{sec:diagonal-sums}

Throughout this section, $p$ denotes a prime occurring in an Euler product; it
is no longer fixed once and for all.

We now pass from a single abelian $p$-group to an arithmetic function of one
integer variable. The relative shape of the invariant factors is fixed, while
the base integer varies. Thus if
\[
        \lambda=(\lambda_1,\ldots,\lambda_{\ell(\lambda)})
\]
is fixed, then, for $n\geq 1$, we consider the groups
\[
        C_{n^{\lambda_1}}\times\cdots\times C_{n^{\lambda_{\ell(\lambda)}}}.
\]
We call this the diagonal problem because the same integer $n$ is used in each
cyclic factor, while the exponents
$\lambda_1,\ldots,\lambda_{\ell(\lambda)}$ remain fixed.

For a fixed nonempty partition $\lambda$, define
\[
        S_\lambda(n)\coloneqq
        s(C_{n^{\lambda_1}}\times\cdots\times C_{n^{\lambda_{\ell(\lambda)}}}).
\]
The function $S_\lambda$ is multiplicative. Indeed, if $(m,n)=1$, then the
Chinese remainder theorem gives
$C_{(mn)^a}\cong C_{m^a}\times C_{n^a}$ for every $a\geq 1$, whence the group
counted by $S_\lambda(mn)$ is a direct product $G_m\times G_n$ of groups of
coprime orders, and every subgroup of such a product splits over the two
factors. Hence
\[
        S_\lambda(mn)=S_\lambda(m)S_\lambda(n).
\]

If $p^b$ is the exact power of $p$ dividing $n$, then the Sylow $p$-subgroup
of $C_{n^{\lambda_i}}$ is $C_{p^{b\lambda_i}}$. Consequently the Sylow
$p$-subgroup of
\[
        C_{n^{\lambda_1}}\times\cdots\times C_{n^{\lambda_{\ell(\lambda)}}}
\]
has type
\[
        b\lambda=(b\lambda_1,\ldots,b\lambda_{\ell(\lambda)}).
\]
Thus
\[
        S_\lambda(p^b)=s(A_{b\lambda}),
\]
where $A_{b\lambda}$ denotes the abelian group of type $b\lambda$ formed over
the prime $p$ under consideration, and $0\lambda$ is interpreted as the empty
partition.

We record the effect of passing from $\lambda$ to $b\lambda$ for $b\geq 1$. If
\[
        \lambda'=(\lambda'_1,\lambda'_2,\ldots,\lambda'_m),
\]
then multiplying the parts of $\lambda$ by $b$ replaces each column of
$\lambda$ by $b$ equal columns. Hence
\[
        (b\lambda)'
        =
        (\underbrace{\lambda'_1,\ldots,\lambda'_1}_{b\textnormal{ times}},
         \underbrace{\lambda'_2,\ldots,\lambda'_2}_{b\textnormal{ times}},
         \ldots,
         \underbrace{\lambda'_m,\ldots,\lambda'_m}_{b\textnormal{ times}}).
\]
It follows that
\[
        \eta(b\lambda)=b\eta(\lambda),
        \qquad
        r(b\lambda)=r(\lambda),
\]
and, for $h\geq 1$,
\[
        d_h(b\lambda)=b\,d_h(\lambda).
\]
Therefore
\[
        \theta(b\lambda)=
        \prod_{\substack{h\geq 1\\ h\textnormal{ odd}}}
        (b\,d_h(\lambda)+1).
\]
Thus $\theta(b\lambda)$ grows at most polynomially in $b$, whereas the exponent
of $p$ is exactly $b\eta(\lambda)$.

\begin{theorem}\label{thm:fixed-shape}
Let $\lambda$ be a fixed nonempty partition, and put
\[
        \eta=\eta(\lambda),
        \qquad
        \theta=\theta(\lambda).
\]
There is a function $H_\lambda(s)$, holomorphic in the half-plane
$\operatorname{Re}(s)>\eta+1/2$, such that, for $\operatorname{Re}(s)>\eta+1$,
\[
        \sum_{n=1}^{\infty}\frac{S_\lambda(n)}{n^s}
        =
        \zeta(s-\eta)^\theta H_\lambda(s).
\]
Thus the right-hand side gives the meromorphic continuation of this Dirichlet
series to $\operatorname{Re}(s)>\eta+1/2$.
More explicitly, in this half-plane,
\[
        H_\lambda(s)\coloneqq
        \prod_p
        (1-p^{\eta-s})^\theta
        \left(\sum_{b=0}^{\infty}\frac{s(A_{b\lambda})}{p^{bs}}\right),
\]
with normal convergence on compact subsets.
Consequently there is a polynomial $P_\lambda$ of degree $\theta-1$ such
that
\[
        \sum_{n\leq x}S_\lambda(n)
        =
        x^{\eta+1}P_\lambda(\log x)
        +
        O_\lambda\left(\frac{x^{\eta+1}}{\log x}\right).
\]
The leading coefficient of $P_\lambda$ is
\[
        \frac{H_\lambda(\eta+1)}{(\eta+1)\Gamma(\theta)},
\]
where $\Gamma$ denotes Euler's Gamma function.
Equivalently,
\[
        \sum_{n\leq x}S_\lambda(n)
        \sim
        \frac{H_\lambda(\eta+1)}{(\eta+1)\Gamma(\theta)}
        x^{\eta+1}(\log x)^{\theta-1}.
\]
\end{theorem}

\begin{proofof}
For $\operatorname{Re}(s)$ sufficiently large, multiplicativity gives the
Euler product
\[
        \sum_{n=1}^{\infty}\frac{S_\lambda(n)}{n^s}
        =
        \prod_p
        \left(\sum_{b=0}^{\infty}\frac{s(A_{b\lambda})}{p^{bs}}\right).
\]
We estimate a local factor. In view of Theorem~\ref{thm:main} and
Corollary~\ref{cor:uniform}, there is a constant $C_\lambda$, independent of
$p$ and $b$, such that, for $b\geq 1$,
\[
        s(A_{b\lambda})
        \leq
        C_\lambda (b+1)^{|\lambda|}p^{b\eta}.
\]
Indeed, by the scaling relations recorded above, $\eta(b\lambda)=b\eta$ and
$r(b\lambda)=r(\lambda)$, while the product formula for $\theta(b\lambda)$
recorded there gives $\theta(b\lambda)\leq C'_\lambda(b+1)^{|\lambda|}$.
Finally, $K_p\leq K_2$ for every prime $p$, because
$(p^{-1};p^{-1})_\infty\geq (2^{-1};2^{-1})_\infty$. Thus the factor
$K_p^{3r(\lambda)}$ in Corollary~\ref{cor:uniform} is bounded by a constant
that depends only on $\lambda$.

Set
\[
        u_p\coloneqq p^{\eta-s}.
\]
If $\operatorname{Re}(s)>\eta+1/2$, then $|u_p|<p^{-1/2}$. The tail beginning
with $b=2$ is locally uniformly small. Indeed, on any substrip
$\operatorname{Re}(s)\geq \eta+1/2+\delta$, with $\delta>0$, we have
$|u_p|\leq 2^{-1/2-\delta}<1$, and therefore
\[
\begin{aligned}
        \sum_{b=2}^{\infty}\frac{s(A_{b\lambda})}{p^{bs}}
        &\ll_\lambda
        \sum_{b=2}^{\infty}(b+1)^{|\lambda|}|u_p|^b                 \\
        &\leq
        |u_p|^2\sum_{m=0}^{\infty}(m+3)^{|\lambda|}2^{-m(1/2+\delta)}
        =O_{\lambda,\delta}(|u_p|^2).
\end{aligned}
\]
Consequently, on each compact subset of the half-plane
$\operatorname{Re}(s)>\eta+1/2$,
\[
        \sum_{b=2}^{\infty}\frac{s(A_{b\lambda})}{p^{bs}}
        =
        O_\lambda(|u_p|^2),
\]
where the implied constant may depend on the compact subset.
For the term $b=1$, the asymptotic assertion in Theorem~\ref{thm:main} gives
\[
        s(A_\lambda)=\theta p^\eta+O_\lambda(p^{\eta-1}).
\]
Therefore
\[
        \frac{s(A_\lambda)}{p^s}
        =
        \theta u_p+O_\lambda(p^{-1}|u_p|).
\]
Combining the terms $b=0$, $b=1$, and $b\geq 2$ gives
\[
        \sum_{b=0}^{\infty}\frac{s(A_{b\lambda})}{p^{bs}}
        =
        1+\theta u_p+O_\lambda(p^{-1}|u_p|+|u_p|^2).
\]
Since $\theta$ is a fixed positive integer,
\[
        (1-u_p)^\theta
        =
        1-\theta u_p+O_\lambda(|u_p|^2).
\]
Thus the $p$-th Euler factor of $H_\lambda(s)$ has the form
\[
        1+O_\lambda(p^{-1}|u_p|+|u_p|^2).
\]
Now
\[
        p^{-1}|u_p|=p^{-1+\eta-\operatorname{Re}(s)}
\]
and
\[
        |u_p|^2=p^{2\eta-2\operatorname{Re}(s)}.
\]
Both sums over primes converge locally uniformly in the half-plane
$\operatorname{Re}(s)>\eta+1/2$, whence the Euler product defining
$H_\lambda(s)$ converges normally in that half-plane, and $H_\lambda$ is
holomorphic there. The same majorants are independent of
$\operatorname{Im}(s)$ and show that, for every $\varepsilon>0$, the product
defining $H_\lambda(s)$ is bounded in the half-plane
$\operatorname{Re}(s)\geq \eta+1/2+\varepsilon$.

At the real point $s=\eta+1$ the quantity $u_p$ equals $p^{-1}$, so the $p$-th
factor of the Euler product for $H_\lambda$, namely
\[
        (1-p^{-1})^\theta
        \sum_{b=0}^{\infty}\frac{s(A_{b\lambda})}{p^{b(\eta+1)}},
\]
is positive. Since the Euler product converges normally at $s=\eta+1$ and its
local factors are positive, it follows that $H_\lambda(\eta+1)>0$; in
particular the leading coefficient recorded below is nonzero, and $P_\lambda$
has degree exactly $\theta-1$.

It remains to pass from the Dirichlet series to a summatory estimate. We apply
the Selberg--Delange theorem, in the form stated in
Tenenbaum~\cite[Chap.~II.5, Th.~5.2]{Tenenbaum}, to the shifted series
\[
        \sum_{n=1}^{\infty}\frac{S_\lambda(n)n^{-\eta}}{n^s}
        =
        \zeta(s)^\theta H_\lambda(s+\eta).
\]
The function $H_\lambda(s+\eta)$ is holomorphic in
$\operatorname{Re}(s)>1/2$. After reducing the constant $c_0$ if necessary, the
region
\[
        \operatorname{Re}(s)\geq
        1-\frac{c_0}{\log(2+|\operatorname{Im}(s)|)}
\]
lies in a half-plane $\operatorname{Re}(s)>1/2+\varepsilon$. The preceding
normal convergence gives boundedness there. Selberg--Delange gives a polynomial
$Q_\lambda$ of degree $\theta-1$ such that
\[
        \sum_{n\leq x}S_\lambda(n)n^{-\eta}
        =
        xQ_\lambda(\log x)+O_\lambda\left(\frac{x}{\log x}\right).
\]
We remove the weight $n^{-\eta}$ by partial summation. Let
\[
        W_\lambda(x)\coloneqq\sum_{n\leq x}S_\lambda(n)n^{-\eta}.
\]
Then
\[
        \sum_{n\leq x}S_\lambda(n)
        =
        x^\eta W_\lambda(x)-\eta\int_1^x t^{\eta-1}W_\lambda(t)\,dt.
\]
Substituting the asymptotic expression for $W_\lambda(t)$ gives a main term
$x^{\eta+1}$ times a polynomial in $\log x$, of the same degree
$\theta-1$. The error term is
\[
        O_\lambda\left(x^\eta\frac{x}{\log x}\right)
        +
        O_\lambda\left(\int_2^x t^{\eta-1}\frac{t}{\log t}\,dt\right)
        =
        O_\lambda\left(\frac{x^{\eta+1}}{\log x}\right),
\]
the bounded lower range being absorbed into the error term; when $\eta=0$ the
integral term is absent. Thus
\[
        \sum_{n\leq x}S_\lambda(n)
        =
        x^{\eta+1}P_\lambda(\log x)
        +O_\lambda\left(\frac{x^{\eta+1}}{\log x}\right).
\]
It remains only to record the leading coefficient. 
Since $\zeta(s-\eta)^\theta$ has a pole of order $\theta$ at $s=\eta+1$, 
the leading coefficient of $Q_\lambda$ is
$\kappa_\lambda\coloneqq H_\lambda(\eta+1)/\Gamma(\theta)$. 
The leading part of $W_\lambda(x)$ is $\kappa_\lambda x(\log x)^{\theta-1}$. 
Hence the leading part after partial summation is
\[
        \kappa_\lambda x^{\eta+1}(\log x)^{\theta-1}
        -
        \eta \kappa_\lambda\int_1^x t^\eta(\log t)^{\theta-1}\,dt.
\]
Since
\[
        \int_1^x t^\eta(\log t)^{\theta-1}\,dt
        =
        \frac{x^{\eta+1}}{\eta+1}(\log x)^{\theta-1}
        +O\bigl(x^{\eta+1}(\log x)^{\theta-2}\bigr),
\]
the leading coefficient is
\[
        \kappa_\lambda\left(1-\frac{\eta}{\eta+1}\right)
        =
        \frac{H_\lambda(\eta+1)}{(\eta+1)\Gamma(\theta)}.
\]
This proves the asserted asymptotic formula.
\end{proofof}

\begin{example}
Let $1\leq a\leq b\leq c$, and take $\lambda=(c,b,a)$. Then
$S_\lambda(n)=s(C_{n^a}\times C_{n^b}\times C_{n^c})$.
The conjugate partition is
\[
        \lambda'=(3^a,2^{b-a},1^{c-b}),
\]
where terms with exponent $0$ are omitted. Hence
\[
        \eta(\lambda)
        =
        a\left\lfloor\frac{3^2}{4}\right\rfloor
        +(b-a)\left\lfloor\frac{2^2}{4}\right\rfloor
        +(c-b)\left\lfloor\frac{1^2}{4}\right\rfloor
        =
        a+b.
\]
The odd column heights are $3$ and $1$. Their multiplicities are $a$ and
$c-b$, respectively, so
\[
        \theta(\lambda)=(a+1)(c-b+1).
\]
Also
\[
        r(\lambda)=
        \begin{cases}
        1, & a=b,\\
        2, & a<b.
        \end{cases}
\]
Theorem~\ref{thm:fixed-shape} therefore gives a polynomial $P_{a,b,c}$ of
degree $(a+1)(c-b+1)-1$ such that
\[
        \sum_{n\leq x}
        s(C_{n^a}\times C_{n^b}\times C_{n^c})
        =
        x^{a+b+1}P_{a,b,c}(\log x)
        +
        O_{a,b,c}\left(\frac{x^{a+b+1}}{\log x}\right).
\]
Thus $\eta(\lambda)$ governs the exponent of $x$, while $\theta(\lambda)$
governs the power of $\log x$.
\end{example}

Taking $\lambda=(1^k)$ gives the diagonal problem for $C_n^k$, which is
mentioned in the final remarks of T\'oth's survey~\cite[Sec.~7]{TothSurvey}.
In this case
\[
        S_\lambda(n)=s(C_n^k),
        \qquad
        \eta(\lambda)=\left\lfloor\frac{k^2}{4}\right\rfloor.
\]
There is only one column height, namely $k$. Therefore
\[
        \theta(\lambda)=
        \begin{cases}
        1, & k\textnormal{ even},\\
        2, & k\textnormal{ odd}.
        \end{cases}
\]

\begin{corollary}\label{cor:diagonal}
For fixed $k\geq 1$, set $S_k(n)\coloneqq s(C_n^k)$ and
$\Delta_k\coloneqq\lfloor k^2/4\rfloor$. If $k$ is even, then
\[
        \sum_{n\leq x}S_k(n)
        \sim
        \mathfrak{C}_k x^{\Delta_k+1}.
\]
If $k$ is odd, then
\[
        \sum_{n\leq x}S_k(n)
        \sim
        \mathfrak{C}_k x^{\Delta_k+1}\log x.
\]
In both cases $\mathfrak{C}_k>0$ is the constant obtained from the Euler
product in Theorem~\ref{thm:fixed-shape} with $\lambda=(1^k)$.
\end{corollary}

\begin{proofof}
Apply Theorem~\ref{thm:fixed-shape} with $\lambda=(1^k)$. The computation
above gives $\eta(\lambda)=\Delta_k$. If $k$ is even, then
$\theta(\lambda)=1$, and $P_\lambda$ has degree $0$, giving a main term of the
form $\mathfrak{C}_k x^{\Delta_k+1}$. If $k$ is odd, then
$\theta(\lambda)=2$, and the leading term of the degree-one polynomial
$P_\lambda$ gives
$\mathfrak{C}_k x^{\Delta_k+1}\log x$. The constant $\mathfrak{C}_k$ is
positive because $H_\lambda(\eta+1)>0$.
\end{proofof}

For $k=1$, the group $C_n$ has one subgroup for each divisor of $n$, so
$S_1(n)=\tau(n)$ and Corollary~\ref{cor:diagonal} recovers the main term in the
classical divisor problem. For $k=2$, one obtains the main term
$(5\pi^2/24)x^2$, in agreement with the formula for rank two recorded in T\'oth's
survey.

\section{Leading terms in fixed rank}
\label{sec:fixed-rank}

Section~\ref{sec:delsarte} shows that the degree of $s(A_\lambda)$, as a
polynomial in $p$, is $\eta(\lambda)$, and that the leading coefficient is
$\theta(\lambda)$. The leading term can therefore be read
directly from the conjugate partition $\lambda'$.

The leading term below for rank four is consistent with the leading term proved
by Chew, Chin and Lim~\cite[Cor.~4.4]{ChewChinLim}. The case of rank three
can also be compared with the explicit formulae for rank three of
Hampejs--T\'oth~\cite{HampejsTothRankThree}.

\begin{corollary}\label{cor:fixed-rank-leading-term}
Let $k\geq 1$ and let $1\leq a_1\leq\cdots\leq a_k$. Put
\[
        f(a_1,\ldots,a_k;p)
        \coloneqq
        s(C_{p^{a_1}}\times\cdots\times C_{p^{a_k}}).
\]
Set $a_0=0$. Then $f(a_1,\ldots,a_k;p)$ is a polynomial in $p$ of degree
\[
        \sum_{j=1}^k a_j\left\lfloor\frac{k-j+1}{2}\right\rfloor
\]
and leading coefficient
\[
        \prod_{\substack{1\leq j\leq k\\ 2\mid k-j}}
        (a_j-a_{j-1}+1).
\]
\end{corollary}

\begin{proofof}
The group $C_{p^{a_1}}\times\cdots\times C_{p^{a_k}}$ corresponds to the
partition
\[
        \lambda=(a_k,a_{k-1},\ldots,a_1).
\]
Its conjugate partition is
\[
        \lambda'
        =
        (k^{a_1},(k-1)^{a_2-a_1},\ldots,1^{a_k-a_{k-1}}),
\]
where terms with exponent $0$ are omitted. Therefore
\[
        \eta(\lambda)
        =
        \sum_{j=1}^k
        (a_j-a_{j-1})
        \left\lfloor\frac{(k-j+1)^2}{4}\right\rfloor.
\]
Put $F(h)\coloneqq\lfloor h^2/4\rfloor$. Since $F(0)=0$ and, for $h\geq 1$,
\[
        F(h)-F(h-1)=\left\lfloor\frac{h}{2}\right\rfloor,
\]
summation by parts gives
$\sum_{j=1}^k(a_j-a_{j-1})F(k-j+1)=
\sum_{j=1}^k a_j(F(k-j+1)-F(k-j))$, and hence
\[
\begin{split}
        \eta(\lambda)
        &=
        \sum_{j=1}^k
        a_j
        \left(
        F(k-j+1)-F(k-j)
        \right)  \\
        &=
        \sum_{j=1}^k
        a_j\left\lfloor\frac{k-j+1}{2}\right\rfloor.
\end{split}
\]

It remains to compute $\theta(\lambda)$. The multiplicity of the part
$k-j+1$ in $\lambda'$ is $a_j-a_{j-1}$. Since $\theta(\lambda)$ is the product
of $d_h(\lambda)+1$ over the odd column heights $h$, we obtain
\[
        \theta(\lambda)
        =
        \prod_{\substack{1\leq j\leq k\\ k-j+1\textnormal{ odd}}}
        (a_j-a_{j-1}+1).
\]
The condition $k-j+1$ odd is equivalent to $2\mid k-j$. Hence
\[
        \theta(\lambda)
        =
        \prod_{\substack{1\leq j\leq k\\ 2\mid k-j}}
        (a_j-a_{j-1}+1).
\]
This gives the asserted degree and leading coefficient.
\end{proofof}

For reference, the cases of ranks $3,4,5$ in
Corollary~\ref{cor:fixed-rank-leading-term} are the following:
\begin{center}
\begin{tabular}{@{}lll@{}}
\toprule
Group & $\deg_p f$ & leading coefficient \\
\midrule
$C_{p^a}\times C_{p^b}\times C_{p^c}$
&
$a+b$
&
$(a+1)(c-b+1)$
\\[2pt]
$C_{p^a}\times C_{p^b}\times C_{p^c}\times C_{p^d}$
&
$2a+b+c$
&
$(b-a+1)(d-c+1)$
\\[2pt]
$C_{p^a}\times C_{p^b}\times C_{p^c}\times C_{p^d}\times C_{p^e}$
&
$2a+2b+c+d$
&
$(a+1)(c-b+1)(e-d+1)$
\\
\bottomrule
\end{tabular}
\end{center}

\section{Open problems}
\label{sec:open-problems}

Theorem~\ref{thm:main} and the three applications above leave several natural
questions open. We single out three that we find appealing, and that we have
not been able to settle.

\begin{question}\label{q:lower-bound}
The lower bound of Theorem~\ref{thm:main} records only the leading term
$\theta(\lambda)p^{\eta(\lambda)}$, whereas the upper bound carries the
homocyclic block constants $\xi_h(q)$ and so responds to the height of each
column block. Is there a lower bound for $s(A_\lambda)$ of comparable
precision, sensitive to the shape of $\lambda$ and equipped with explicit
block constants of its own?
\end{question}

\begin{question}\label{q:small-prime}
The comparison theorem of Section~\ref{sec:class-less-than-p} is, by its very
hypothesis, silent at the prime $2$, since a $2$-group of nilpotency class
less than $2$ is already abelian. Likewise, for $p$-groups of exponent $p$ it
gives only the elementary abelian comparison group, and hence no shape-sensitive
refinement. Can one obtain meaningful subgroup bounds for these two classes?
\end{question}

\begin{question}\label{q:arbitrary-class}
Groves's reduction requires the nilpotency class to be less than $p$. Is
there anything useful to be said, in the spirit of Theorem~\ref{thm:main},
about the subgroup count of a $p$-group of arbitrary nilpotency class?
\end{question}

We would be glad to see progress on any of the three.

\section*{Acknowledgements}

It is a pleasure to thank Benjamin Sambale, who drew our attention to
Groves's paper~\cite{Groves} as well as T\'oth's survey~\cite{TothSurvey} in a private communication.

\end{document}